\documentclass[11pt,a4paper]{article}
\usepackage[T1]{fontenc}
\usepackage[utf8]{inputenc}
\usepackage{lmodern}
\usepackage[margin=25mm]{geometry}
\usepackage{amsmath,amssymb,amsfonts,amsthm}
\usepackage{graphicx,enumitem}
\usepackage{tikz}
\usetikzlibrary{arrows.meta}
\usepackage{float,placeins,needspace}
\usepackage[numbers,sort&compress]{natbib}
\usepackage[font=small,labelfont=bf]{caption}
\usepackage{microtype}
\usepackage[hidelinks,hypertexnames=false]{hyperref}

\definecolor{figorange}{HTML}{EB6834}
\definecolor{figink}{HTML}{0B0B0B}
\definecolor{figinktwo}{HTML}{52514E}
\definecolor{figpaleblue}{HTML}{DCE8F7}
\definecolor{figpaleboth}{HTML}{F1ECE6}

\theoremstyle{plain}
\newtheorem{theorem}{Theorem}[section]
\newtheorem{lemma}[theorem]{Lemma}
\newtheorem{proposition}[theorem]{Proposition}
\newtheorem{corollary}[theorem]{Corollary}
\theoremstyle{definition}
\newtheorem{definition}[theorem]{Definition}
\newtheorem{remark}[theorem]{Remark}
\newtheorem{example}[theorem]{Example}
\numberwithin{equation}{section}

\setlist{itemsep=2pt,topsep=5pt}
\hypersetup{
  pdftitle={Iterated minimal Darboux transformations},
  pdfauthor={Sakuma Takeshita and Masashi Yasumoto},
  bookmarksdepth=2,
  allcolors=black
}
\title{Iterated minimal Darboux transformations}
\author{Sakuma Takeshita \and Masashi Yasumoto}
\date{}

\begin{document}
\maketitle

\begin{abstract}
In this paper, we investigate the superposition of two minimal Darboux transformations in Euclidean three-space. A four-dimensional permutability diagram for transformations of isothermic surfaces yields the common minimal transform for distinct nonzero parameters and an explicit formula for its Gauss map. This formula determines the planar ends over regular umbilic-free points of the seed: an end is embedded when the two transformed Gauss maps coincide away from the seed value and has multiplicity 3 when all three values coincide. We illustrate this criterion with the Enneper surface and the catenoid. For the catenoid, suitable parameters give complete surfaces of finite total curvature, and we determine the numbers and multiplicities of their planar ends from the parameters and initial values. In the family with multiplicity-3 ends, the remaining initial-value parameter realizes, up to rotation, the L\'opez--Ros deformation.
\end{abstract}



\allowdisplaybreaks[0]

\section{Introduction}
Transformations and deformations are classical tools for constructing minimal surfaces in Euclidean three-space $\mathbb{R}^3$ from known ones. The associated family, the Goursat transformation~\citep{Goursat} and the L\'opez--Ros deformation~\citep{LR} act on the Gauss map and holomorphic one-form in the Weierstrass representation~\citep{Osserman}. They reveal relations between seemingly different surfaces and also serve as tools in the global theory; see, for example, \citep{LR,LM,AkamineFujino}.

The transformation theory of isothermic surfaces provides another approach. Isothermic surfaces admit conformal curvature line coordinates. Away from umbilics, examples include minimal surfaces and surfaces with constant nonzero mean curvature. Their classical transformations include those of Christoffel~\citep{Christoffel}, Goursat~\citep{Goursat,HJ} and Darboux~\citep{Darboux}. Bianchi's permutability theorem constructs the common transform of two Darboux transforms algebraically, without further integration~\citep{Bianchi,Bianchi1905}; see also~\citep{HJ,HP}. Darboux transformations can produce substantial changes in geometry: for example, transforms of the round cylinder~\citep{HP,SW} include constant mean curvature bubbletons. Cho, Leschke and Ogata~\citep{CLO} extended Bianchi permutability to equal parameters and used it to construct further bubbletons.

A Darboux transform of a minimal surface need not be minimal. Transformations that preserve minimality, called \emph{minimal Darboux transformations}, were constructed as special Ribaucour transformations by Corro, Ferreira and Tenenblat~\citep{CFT}. Mart\'inez, Roitman and Tenenblat~\citep{MRT} characterized them by a Riccati equation for the Gauss maps. Hertrich-Jeromin and Honda~\citep{HH} recovered this characterization through the permutability of Christoffel, Goursat, Darboux and M\"obius transformations. Minimal Darboux transformations can create planar ends, making their iteration a way to construct minimal surfaces with additional ends. Raujouan, Schmitt and Ziefle~\citep{RSZ} obtained such transformations of the catenoid by simple factor dressing in the loop Weierstrass representation. Here we investigate how two minimal Darboux transformations can be superposed and how their initial values determine the resulting ends.

Starting from the classical Bianchi theorem, we combine two Darboux directions with Christoffel, Goursat and M\"obius transformations in a four-dimensional permutability scheme (Theorem~\ref{p1:diagram}). For distinct nonzero parameters, this scheme shows that the common Darboux transform is again minimal on its regular set (Theorem~\ref{p1:minimal}). It also gives an explicit algebraic formula for its Gauss map in terms of the Gauss maps of the seed and its two initial transforms; see~\eqref{p1:superposition}. Thus another Riccati equation need not be solved. Furthermore, the construction recovers the superposition formula for solutions of the Liouville equation given in~\citep[Theorem~3.3]{CFT}, with a formula valid for distinct nonzero real parameters of either sign (Corollary~\ref{p1:liouville}).

The explicit Gauss map allows us to determine the end multiplicities. Mart\'inez, Roitman and Tenenblat~\citep{MRT} analyzed the effects of a single transformation on regular points, umbilics and ends, including the occurrence of non-embedded planar ends. For the superposition, we show that an end lies over a regular umbilic-free point of the seed exactly when the two transformed Gauss maps coincide (Proposition~\ref{p1:branchorders} and Corollary~\ref{p1:endcor}). It is embedded and planar if their common value differs from the seed value and has multiplicity~3 if all three values coincide. In the latter case, both intermediate surfaces have embedded planar ends there. Our calculation therefore identifies the initial values producing a non-embedded end and determines its multiplicity in the superposition.

We illustrate these results with the Enneper surface and the catenoid. For the catenoid, suitable integer values of $\lambda_1,\lambda_2$ give complete surfaces of finite total curvature with two catenoidal ends and finitely many planar ends. Proposition~\ref{p1:global} determines the numbers and multiplicities of the planar ends, and hence the total curvature, from the parameters and initial values; Figure~\ref{p1:end-count-figure} summarizes the alternatives. Under condition~\eqref{p1:cat-normalization}, the remaining initial-value parameter realizes, up to rotation, the L\'opez--Ros deformation of the superposed surface (Proposition~\ref{p1:lopezros}). Thus this classical deformation appears within the superposition family through the choice of initial values.
\begin{figure}[H]

\centering
\resizebox{\linewidth}{!}{%
\begin{tikzpicture}[x=4.5cm, y=4.5cm,
  lattice/.style={figink, line width=0.45pt, dash pattern=on 0.45pt off 1pt},
  axis/.style={-{Latex[length=5pt,width=4pt]}, figink, line width=0.8pt},
  orangeline/.style={figorange, line width=2.2pt, line cap=round},
  orangedot/.style={circle, fill=figorange, draw=white, line width=0.8pt, inner sep=0pt, minimum size=5.5pt},
  leader/.style={figinktwo, line width=0.5pt},
  leaderdot/.style={circle, fill=figinktwo, inner sep=0pt, minimum size=2.4pt},
  tick/.style={figink, line width=0.6pt}]
  \def\m{0.25}

  \begin{scope}
    \fill[figpaleblue] (\m,\m) rectangle (1,1);
    \foreach \k in {3,...,7} {
      \draw[lattice] (\k/8,\m) -- (\k/8,1);
      \draw[lattice] (\m,\k/8) -- (1,\k/8);
    }
    \draw[figinktwo, line width=0.7pt] (\m,\m) rectangle (1,1);
    \draw[axis] (0,0) -- (1.07,0) node[right, inner sep=2pt] {$\lambda_1$};
    \draw[axis] (0,0) -- (0,1.09) node[above, inner sep=2pt] {$\lambda_2$};
    \draw[tick] (\m,0) -- (\m,-0.025) node[below, inner sep=1.5pt, font=\footnotesize] {$2$};
    \draw[tick] (0,\m) -- (-0.025,\m) node[left, inner sep=1.5pt, font=\footnotesize] {$2$};

    \draw[figinktwo, line width=1pt, dash pattern=on 4pt off 2pt] (\m,\m) -- (1,1);
    \node[anchor=south east, inner sep=0pt, text=figinktwo, font=\footnotesize]
      at (1,1.025) {$\lambda_1=\lambda_2$ (excluded)};

    \node[fill=white, rounded corners=3pt, align=center, font=\small]
      at (0.625,0.625) {only embedded\\ planar ends:\\ $\lambda_1+\lambda_2$};
    \node[anchor=west, inner sep=0pt, font=\small] at (0,-0.35)
      {(a) condition (4.7) fails};
  \end{scope}

  \begin{scope}[shift={(1.275,0)}]
    \fill[figpaleboth] (\m,\m) rectangle (1,1);
    \foreach \k in {3,...,7} {
      \draw[lattice] (\k/8,\m) -- (\k/8,1);
      \draw[lattice] (\m,\k/8) -- (1,\k/8);
    }
    \draw[figinktwo, line width=0.7pt] (\m,\m) rectangle (1,1);
    \draw[axis] (0,0) -- (1.07,0) node[right, inner sep=2pt] {$\lambda_1$};
    \draw[axis] (0,0) -- (0,1.09) node[above, inner sep=2pt] {$\lambda_2$};
    \draw[tick] (\m,0) -- (\m,-0.025) node[below, inner sep=1.5pt, font=\footnotesize] {$2$};
    \draw[tick] (0,\m) -- (-0.025,\m) node[left, inner sep=1.5pt, font=\footnotesize] {$2$};

    \draw[figinktwo, line width=1pt, dash pattern=on 4pt off 2pt] (\m,\m) -- (1,1);
    \node[anchor=south west, inner sep=0pt, text=figinktwo, font=\footnotesize]
      at (1.03,1.01) {$\lambda_1=\lambda_2$ (excluded)};

    \coordinate (src) at (1.085,0.285);
    \foreach \p in {(0.85,0.32),(0.87,0.62),(0.62,0.87),(0.31,0.94)} {
      \draw[leader] (src) -- \p;
      \node[leaderdot] at \p {};
    }

    \draw[orangeline] (2/8,4/8) -- (0.5,1);
    \draw[orangeline] (4/8,2/8) -- (1,0.5);
    \node[anchor=south, inner sep=0pt, font=\small] at (0.5,1.035) {$\lambda_2=2\lambda_1$};
    \node[anchor=west, inner sep=0pt, font=\small] at (1.04,0.5) {$\lambda_1=2\lambda_2$};
    \foreach \k in {2,...,4} {
      \node[orangedot] at (\k/8,2*\k/8) {};
      \node[orangedot] at (2*\k/8,\k/8) {};
    }

    \node[anchor=west, align=left, inner sep=0pt, font=\footnotesize] at (1.11,0.2)
      {both:\\ $\lambda_1+\lambda_2-3d$ embedded,\\ $d$ of multiplicity $3$};

    \draw[orangeline] (0,-0.135) -- (0.2,-0.135);
    \node[anchor=west, inner sep=0pt, font=\footnotesize] at (0.26,-0.135)
      {only planar ends of multiplicity $3$:\ \ $d$};
    \node[anchor=west, inner sep=0pt, font=\footnotesize] at (0.26,-0.24)
      {$d=\gcd(\lambda_1,\lambda_2)$};

    \node[anchor=west, inner sep=0pt, font=\small] at (0,-0.35)
      {(b) condition (4.7) holds};
  \end{scope}

\end{tikzpicture}
}%

 \caption{Types and numbers of planar ends of the superposed catenoid $f^{12}$ (Proposition~\ref{p1:global}). The distinct integers $\lambda_1,\lambda_2\geq2$ determine the Darboux parameters. In~(a), condition~\eqref{p1:cat-normalization} holds for no $\alpha\in\mathbb C^*$; in~(b), it holds for some $\alpha\in\mathbb C^*$. Here $d=\gcd(\lambda_1,\lambda_2)$.}
\label{p1:end-count-figure}
\end{figure}

This paper is organized as follows. Section~\ref{sec:preliminaries} reviews isothermic surfaces and their transformations. Section~\ref{sec:superposition} proves the superposition principle, analyzes the ends of the superposed surface, and compares the resulting metric formula with that of~\citep{CFT}. Section~\ref{sec:examples} develops the examples and their global properties.

\section{Preliminaries}\label{sec:preliminaries}
In this section, we use quaternionic notation to recall the definitions of isothermic surfaces and their transformations, together with the permutability relations used below; see~\citep{HJ} for details. Following~\citep{HH}, we then characterize minimal Darboux transformations by a Riccati equation for their Gauss maps.

We first identify $\mathbb R^3$ with the imaginary quaternions. Let $\mathbb H$ be the quaternion algebra. Then $\mathbb R^3$ is identified with
\[
\operatorname{Im}\mathbb H=\{x_1\mathbf i+x_2\mathbf j+x_3\mathbf k:x_1,x_2,x_3\in\mathbb R\},
\qquad \mathbf i^2=\mathbf j^2=\mathbf k^2=\mathbf i\mathbf j\mathbf k=-1.
\]
For $q=a_0+a_1\mathbf i+a_2\mathbf j+a_3\mathbf k$, set
\[
\overline q=a_0-a_1\mathbf i-a_2\mathbf j-a_3\mathbf k,\quad
\operatorname{Re}q=a_0,\quad
\operatorname{Im}q=q-a_0,\quad
|q|^2=q\overline q=\sum_{\ell=0}^3a_\ell^2.
\]
Then $q^{-1}=\overline q/|q|^2$ for $q\ne0$; in particular, $x^{-1}=-x/|x|^2$ for $x\in\operatorname{Im}\mathbb H\setminus\{0\}$. Since
\[
x_1\mathbf i+x_2\mathbf j+x_3\mathbf k=x_1\mathbf i+(x_2+x_3\mathbf i)\mathbf j,
\]
we also identify $\operatorname{Im}\mathbb H$ with $\mathbb R\mathbf i\oplus\mathbb C\mathbf j$, where the quaternionic copy of $\mathbb C$ is $\mathbb R\oplus\mathbb R\mathbf i$.
The symbol ``$\mathrm i$'' denotes the complex imaginary unit in scalar formulas. This identification relates a spherical Gauss map to its stereographic coordinate.

For four points $A,B,C,D\in\mathbb R^3\cong\operatorname{Im}\mathbb H$, the quaternionic cross ratio is
\[
\mathcal Q(A,B,C,D)=(A-B)(B-C)^{-1}(C-D)(D-A)^{-1},
\]
and we represent its conjugacy class by
\[
\operatorname{cr}(A,B,C,D)=\operatorname{Re}\mathcal Q+\mathrm i\,|\operatorname{Im}\mathcal Q|.
\]
For four distinct points, this value is real exactly when the points are concircular or collinear. It is M\"obius invariant, and for points in a complex plane it agrees with the classical cross ratio up to complex conjugation. For scalar complex values we use
\[
\operatorname{cr}(a,b,c,d)=\frac{(a-b)(c-d)}{(b-c)(d-a)};
\]
only real cross ratios are used below, so the conventions agree there.

For vectors $x,y\in\mathbb R^3$, regarded as imaginary quaternions,
\[
xy=-\langle x,y\rangle+x\times y.
\]
M\"obius transformations of $\mathbb R^3\cup\{\infty\}$ act on $\operatorname{Im}\mathbb H\cup\{\infty\}$ by quaternionic fractional linear maps~\citep{HJ}; we only use that they preserve real cross ratios.

\subsection{Isothermic surfaces and their transformations}
Here we define isothermic surfaces in $\mathbb{R}^3$ and three types of transformations of these surfaces. Throughout this paper, Darboux parameters are defined relative to a fixed quadratic differential, locally $dz^2$, which we call the polarization.

\begin{definition}\label{p1:isothermic}
Let $\Sigma$ be a Riemann surface and $f:\Sigma\to\mathbb R^3$ an immersion with unit normal $\nu$, first fundamental form $I_f$, second fundamental form $II_f$, and principal curvatures $\kappa_1,\,\kappa_2$. If there are local conformal curvature line coordinates $(u,v)$ and a real-valued function $\omega$ with
\begin{equation}
I_f=\langle df,df\rangle=e^{2\omega}(du^2+dv^2),\qquad
II_f=-\langle df,d\nu\rangle=e^{2\omega}(\kappa_1du^2+\kappa_2dv^2),
\label{p1:isothermic-metric}
\end{equation}
then $f$ is called an \emph{isothermic surface}.
\end{definition}

It is known that, away from umbilics, minimal surfaces and constant mean curvature surfaces are isothermic. We assume that the minimal seeds are umbilic-free on their regular domains. Auxiliary isothermic surfaces, including the spherical Christoffel duals, may have umbilics; their coordinates are inherited from the chosen polarization. All transformation formulas initially hold on common regular domains with nonzero denominators and compatible integration constants. We subsequently extend the Gauss data meromorphically.

Christoffel~\citep{Christoffel} showed that an isothermic surface $f$ locally admits a surface $f^*$ with the following properties (see also~\citep{HJ}):
\begin{itemize}
\item $f$ and $f^*$ have parallel tangent planes at corresponding points;
\item their shape operators commute, so their curvature lines correspond;
\item their metrics are conformally equivalent, with the induced map reversing orientation.
\end{itemize}
Such a surface $f^*$ is called a Christoffel transform of $f$. With the conformal curvature line coordinates fixed, it is given, in quaternionic form~\citep{HJ}, by
\begin{equation}
f^*_u=(f_u)^{-1},\qquad f^*_v=-(f_v)^{-1}.
\label{p1:christoffel}
\end{equation}
For a fixed polarization it is unique up to translation; allowing a rescaling of the polarization also allows a homothety. It is isothermic, and $(f^*)^*=f$ with suitable normalizations.

Let $z=u+\mathrm i v$. The Hopf differential of $f$ is
\begin{equation}
Q:=\langle f_{zz},\nu\rangle\,dz^2=:q\,dz^2.
\label{p1:hopf}
\end{equation}
For a minimal surface, $q$ is a nonzero real constant on an umbilic-free curvature line patch, and after a rescaling of $z$ we may assume that $q=1/2$; see~\citep{HJ}.

\begin{proposition}[\citep{Christoffel,HJ}]\label{p1:minimaldual}
An umbilic-free isothermic surface is minimal if and only if, after a homothety and translation, its Christoffel transform takes values in $S^2$. In this case $f^*=\nu$ with the normalization $q=1/2$.
\end{proposition}

The classical Goursat transformation of a minimal surface applies a complex orthogonal transformation $A\in O(3,\mathbb C)$ to its holomorphic null curve and acts on its Gauss map by a M\"obius transformation of $\mathbb C\cup\{\infty\}$~\citep{Goursat}. Hertrich-Jeromin extended this transformation to isothermic surfaces~\citep{HJ}.

\begin{definition}\label{p1:goursat}
Let $f$ be isothermic, $f^*$ a Christoffel transform, and $\mu$ a M\"obius transformation. The surface
\[
\tilde f:=(\mu\circ f^*)^*
\]
is called a Goursat transform and is denoted by $G_\mu(f)$.
\end{definition}

The Goursat transformation yields the Weierstrass representation; see \citep{HH, HJ} for more details. Let
\[
\iota(x):=\mathbf i-2(x-\mathbf i)^{-1}.
\]
This is an involutive M\"obius transformation; it restricts to the stereographic projection $S^2\to\mathbb C\mathbf j\cup\{\infty\}$ from $\mathbf i$ and to its inverse on $\mathbb C\mathbf j$. For a locally unbranched holomorphic function $g$, the Goursat transform
\[
f:=G_\iota\!\left(-\mathbf j\int\frac{dz}{g_z}\right)
\]
is a minimal surface with Gauss map $\nu=f^*=\iota(g\mathbf j)$~\citep{HH}. We identify $\operatorname{Im}\mathbb H$ with $\mathbb R^3$ by
\[
x_1\mathbf i+x_2\mathbf j+x_3\mathbf k\longmapsto(x_2,x_3,x_1),
\qquad\text{so that}\qquad
\nu=\frac{(2\operatorname{Re}g,2\operatorname{Im}g,|g|^2-1)}{1+|g|^2}.
\]
Then, up to translation, $f$ is given by the classical Weierstrass representation
\begin{equation}
f=\operatorname{Re}\int^z\frac12(1-g^2,\mathrm i(1+g^2),2g)^T h\,dz,
\qquad h=-\frac1{g_z}.
\label{p1:weierstrass}
\end{equation}
Thus the classical Weierstrass data are $(g,h\,dz)=(g,-dz/g_z)$~\citep{Osserman}. In the coordinate normalized by $q=1/2$, the Gauss map determines $f$ up to translation, subject to the period condition if the domain is not simply connected. The fundamental forms are
\begin{equation}
I_f=\frac{(1+|g|^2)^2}{4|g_z|^2}\,dz\,d\overline z,
\qquad II_f=\frac{dz^2+d\overline z^{\,2}}2=du^2-dv^2.
\label{p1:fundamental}
\end{equation}
In particular, $Q=dz^2/2$.

Following Darboux~\citep{Darboux}, we call two surfaces $f,\hat f$ a nontrivial Darboux pair if they envelop a common tangent sphere congruence, have corresponding curvature lines, and induce conformally equivalent metrics. On the nontrivial envelope branch, an isothermic surface $f$ and its Darboux transform $\hat f$ satisfy the quaternionic Riccati equation~\citep{HJ,HP}
\begin{equation}
d\hat f=t(\hat f-f)\,df^*\,(\hat f-f),
\qquad t\in\mathbb R\setminus\{0\}.
\label{p1:qriccati}
\end{equation}
For fixed $t$, the initial value $\hat f(p_0)$ determines a local transform; we write $D_t(f)$ for a selected member of this three-parameter family. Darboux transformations commute with M\"obius transformations~\citep{HJ}, with the induced polarization, whereas Christoffel transformations do not.

\subsection{Permutability of transformations for isothermic surfaces}
In this setting, the following permutability theorem holds~\citep{Bianchi,Bianchi1905,HP}.
\begin{proposition}[Bianchi permutability]\label{p1:bianchiprop}
Let $f$ be an isothermic surface, and let $f^i=D_{t_i}(f) \ (i=1,2)$ be its Darboux transforms with distinct nonzero real parameters $t_i$. On the nondegenerate common domain, there is a unique fourth surface $f^{12}$ that is a Darboux transform of both $f^1$ and $f^2$:
\[
f^{12}=D_{t_1}(f^2)=D_{t_2}(f^1).
\]
It is given explicitly by
\begin{equation}
(f^{12}-f)^{-1}=\frac{t_2(f^1-f)^{-1}-t_1(f^2-f)^{-1}}{t_2-t_1}.
\label{p1:bianchi}
\end{equation}
Furthermore, the fourth surface $f^{12}$ satisfies
\[
\operatorname{cr}(f,f^1,f^{12},f^2)=\frac{t_2}{t_1},
\]
and the four corresponding points lie on a circle or a straight line.
\end{proposition}
We call the diagram of these four surfaces, in the boundary order $(f,f^1,f^{12},f^2)$, a \emph{Darboux square}. The edges from $f$ to $f^1$ and from $f^2$ to $f^{12}$ carry the parameter $t_1$; the other two carry $t_2$. The cross ratio of the square is evaluated at corresponding points of the four surfaces.

The Christoffel transform of a Darboux transform and the Darboux transform of a Goursat transform can be obtained algebraically.

\begin{proposition}\label{p1:cdprop}
Let $f$ be an isothermic surface and let $f^*$ be a Christoffel transform of $f$. Let $\hat f=D_t(f)$ be a Darboux transform of $f$ with parameter $t$, and let $\tilde f=G_\mu(f)$ be a Goursat transform of $f$ for a M\"obius transformation $\mu$. Then, for compatible choices of the integration constants,
\[
\hat f^*=(D_t(f))^*=D_t(f^*),\qquad
\check f=D_t(G_\mu(f))=G_\mu(D_t(f)),
\]
where
\[
\hat f^*=f^*+\frac1t(\hat f-f)^{-1},\qquad
\check f=\tilde f+\frac1t\bigl(\mu(\hat f^*)-\mu(f^*)\bigr)^{-1}.
\]
\end{proposition}

The permutability of Christoffel and Darboux transformations is due to Bianchi~\citep{Bianchi,Bianchi1905}; see also~\citep[Theorem~5.6.3]{HJ} and~\citep[Theorem~2.8]{Burstall}. Hertrich-Jeromin and Honda~\citep{HH} established the extension to Goursat transformations. These identities require compatible relative translations, rather than independently chosen integration constants.

When $f$ and $\hat f$ are minimal and $\mu=\iota$, the corresponding vertices $\iota(f^*)$ and $\iota(\hat f^*)$ are $g\mathbf j$ and $\hat g\mathbf j$. Accordingly, two nonplanar minimal surfaces with Hopf differential $dz^2/2$ in a common coordinate form a Darboux pair, after a suitable relative translation, exactly when their Gauss maps satisfy the Riccati equation~\citep{HH,MRT}
\begin{equation}
\hat g_z=t\frac{(\hat g-g)^2}{g_z},
\qquad t\in\mathbb R\setminus\{0\},
\label{p1:minriccati}
\end{equation}
and, with this choice of translation, their pointwise difference is
\begin{equation}
\hat f-f=\frac1t(\hat\nu-\nu)^{-1},
\qquad \nu=f^*,\quad\hat\nu=\hat f^*.
\label{p1:relative}
\end{equation}
Such a Darboux transformation is called a minimal Darboux transformation.

\section{Superposition principle for minimal Darboux transformations}\label{sec:superposition}
In this section, we combine two Darboux transformations with Christoffel, Goursat and M\"obius transformations in a compatible diagram. This gives the superposition principle for minimal Darboux transformations and the Gauss-map formula used to analyze the ends of the superposed surface. Throughout this paper, $t_1,t_2$ are distinct nonzero real parameters, and $f^{12}$ is the common transform in Proposition~\ref{p1:bianchiprop} on its regular domain. Let $f$ be isothermic, $f^*$ its Christoffel transform, and $f^i=D_{t_i}(f)$. Set
\[
f^{i*}:=f^*+\frac1{t_i}(f^i-f)^{-1},\qquad i=1,2,
\]
as in Proposition~\ref{p1:cdprop}. Bianchi~\citep{Bianchi1905} showed that taking Christoffel duals at the four vertices gives another Darboux square, which we call the \emph{Christoffel-dual square}; see \citep[Theorem 2.20 and Exercise 2.9]{Burstall} for an algebraic proof. Indeed, the surface
\begin{equation}
f^{12*}:=f^{1*}+\frac1{t_2}(f^{12}-f^1)^{-1}
=f^{2*}+\frac1{t_1}(f^{12}-f^2)^{-1}
\label{p1:dualclosure}
\end{equation}
is a Christoffel transform of $f^{12}$, and
\[
f^{12*}=D_{t_2}(f^{1*})=D_{t_1}(f^{2*}),\qquad
\operatorname{cr}(f^*,f^{1*},f^{12*},f^{2*})=\frac{t_2}{t_1}.
\]
For a M\"obius transformation $\mu$, we call the image of the Christoffel-dual square under $\mu$ the \emph{M\"obius square}. Applying $G_\mu$ to the four original vertices gives the \emph{Goursat square}. The next lemma fixes the translation of $G_\mu(f^{12})$.

\begin{lemma}\label{p1:goursatlemma}
Let $\mu$ be a M\"obius transformation, choose $G_\mu(f^i)$ by Proposition~\ref{p1:cdprop} and put
\[
G_\mu(f^{12}):=G_\mu(f^1)+\frac1{t_2}
\bigl(\mu(f^{12*})-\mu(f^{1*})\bigr)^{-1}.
\]
Then
\[
\begin{aligned}
G_\mu(f^{12})
&=G_\mu(f^2)+\frac1{t_1}\bigl(\mu(f^{12*})-\mu(f^{2*})\bigr)^{-1}
=(\mu\circ f^{12*})^*,\\
G_\mu(f^{12})&=D_{t_2}(G_\mu(f^1))=D_{t_1}(G_\mu(f^2)),
\end{aligned}
\]
and the Goursat square has cross ratio $t_2/t_1$.
\end{lemma}
\begin{proof}
Since Darboux transformations and real cross ratios are M\"obius invariant, \eqref{p1:dualclosure} gives
$\mu(f^{12*})=D_{t_2}(\mu(f^{1*}))=D_{t_1}(\mu(f^{2*}))$ and
$\operatorname{cr}(\mu(f^*),\mu(f^{1*}),\mu(f^{12*}),\mu(f^{2*}))=t_2/t_1$;
thus the M\"obius square is a Darboux square.
Its first vertex has the Christoffel transform $G_\mu(f)=(\mu\circ f^*)^*$, and by Proposition~\ref{p1:cdprop} the corresponding Christoffel transforms of $\mu(f^{i*})=D_{t_i}(\mu(f^*))$ are the surfaces $G_\mu(f^i)$. Applying \citep[Theorem 2.20]{Burstall}, as in \eqref{p1:dualclosure}, to the M\"obius square gives the assertions.
\end{proof}

In Figure~\ref{p1:permutability}(b), two directions follow the Darboux edges and the other two follow the horizontal and vertical links between squares. A \emph{face} of the full diagram is obtained by varying two directions and fixing the other two. It commutes if the two paths from its initial vertex to the opposite vertex give the same surface. For instance, the face through $f$ spanned by $D_{t_1}$ and $G_\mu$ records $G_\mu(f^1)=D_{t_1}(G_\mu(f))$.

\begin{theorem}\label{p1:diagram}
In the setting of Lemma~\ref{p1:goursatlemma}, $f^{12*}$ and $G_\mu(f^{12})$ are isothermic, and all faces of the diagram in Figure~\ref{p1:permutability} commute. Its edges represent Christoffel, Goursat, M\"obius, or Darboux transformations. The four Darboux squares in panel~(b) all have cross ratio $t_2/t_1$.
\end{theorem}
\begin{proof}
The top row of Figure~\ref{p1:permutability}(b) consists of the original Darboux square of Proposition~\ref{p1:bianchiprop} and its Christoffel-dual square. Proposition~\ref{p1:cdprop} and~\eqref{p1:dualclosure} show that Christoffel duality preserves both Darboux directions. Applying $\mu$ to the dual square gives the M\"obius square (lower right):
\begin{equation}
\begin{aligned}
\operatorname{cr}(\mu(f^*),\mu(f^{1*}),\mu(f^{12*}),\mu(f^{2*}))&=\frac{t_2}{t_1},\\
\mu(f^{12*})&=D_{t_1}(\mu(f^{2*}))=D_{t_2}(\mu(f^{1*})).
\end{aligned}
\label{p1:mobiusface}
\end{equation}
Lemma~\ref{p1:goursatlemma} gives the Goursat square (lower left), whose edges correspond to those of the original square by Proposition~\ref{p1:cdprop}. Finally, $(G_\mu(f))^*=\mu(f^*)$ with the chosen translations, and the same identity holds at $f^1,f^2,f^{12}$. These identities close the faces spanned by the horizontal and vertical links in panel~(b).
\end{proof}

We now apply Theorem~\ref{p1:diagram} to a minimal surface $f$. The Christoffel-dual square shows that $f^{12}$ is again minimal.
\begin{theorem}\label{p1:minimal}
Let $(f,f^i)$, $i=1,2$, be minimal Darboux pairs with distinct nonzero parameters. Then the common Darboux transform $f^{12}$ of Proposition~\ref{p1:bianchiprop} is minimal on its regular set. Thus $(f^i,f^{12})$ are minimal Darboux pairs there.
\end{theorem}

\begin{remark}
In~\cite{TY}, we give an alternative proof of Theorem~\ref{p1:minimal} using the Weierstrass representation of Bobenko and Pinkall~\cite{BP} for discrete minimal surfaces. 
\end{remark}

\begin{proof}[Proof of Theorem~\ref{p1:minimal}]
With compatible normalizations, $f^*,f^{1*},f^{2*}$ are the Gauss maps of $f,f^1,f^2$~\citep[Corollary 3.3]{HH}, so they lie on $S^2$. By~\eqref{p1:dualclosure}, the Christoffel-dual square has real cross ratio $t_2/t_1$.
The three vertices $f^*,f^{1*},f^{2*}$ are distinct on a dense open set. Equality of $f^{i*}$ with $f^*$ is excluded by~\eqref{p1:relative}. If $f^{1*}=f^{2*}$ on an open set, choose a finite spherical chart there. Then $g^1=g^2$, and~\eqref{p1:minriccati} gives
\[
(t_1-t_2)(g^1-g)^2=0,
\]
contradicting $t_1\ne t_2$ and~\eqref{p1:relative}. Real analyticity now gives the claimed dense open set. The circle through these three vertices is the intersection of $S^2$ with their plane, so $f^{12*}\in S^2$ there, and hence throughout its domain by continuity. Proposition~\ref{p1:minimaldual} gives minimality of $f^{12}$ on its regular set.
\end{proof}

\begin{figure}[tbp]
\centering
\begin{tikzpicture}[
  edge/.style={line width=0.5pt},
  hidden/.style={line width=0.5pt, dash pattern=on 3pt off 3pt},
  vtx/.style={circle, fill, inner sep=0pt, minimum size=2.8pt},
  sq/.style={line width=0.55pt, -{Stealth[length=4pt,width=3pt,inset=1.3pt]}},
  link/.style={double, double distance=0.8pt, line width=0.35pt, -{Implies}},
  lab/.style={inner sep=2pt}]

\begin{scope}
  \coordinate (Gf)   at (0,0);     \coordinate (Mf)   at (4,0);
  \coordinate (Gf2)  at (0,4);     \coordinate (Mf2)  at (4,4);
  \coordinate (Gf1)  at (2,1.2);   \coordinate (Mf1)  at (6,1.2);
  \coordinate (Gf12) at (2,5.2);   \coordinate (Mf12) at (6,5.2);
  \coordinate (f)   at (1.6,1.6);  \coordinate (s)   at (3.6,1.6);
  \coordinate (f2)  at (1.6,3.2);  \coordinate (s2)  at (3.6,3.2);
  \coordinate (f1)  at (2.4,2.24); \coordinate (s1)  at (4.4,2.24);
  \coordinate (f12) at (2.4,3.84); \coordinate (s12) at (4.4,3.84);

  \draw[edge] (f)--(s) (f2)--(s2) (f12)--(s12) (f1)--(s1);
  \draw[edge] (f)--(f1)--(f12);
  \draw[edge] (s)--(s1) (f2)--(f12) (s2)--(s12);
  \draw[edge] (f)--(f2) (s)--(s2) (s1)--(s12);
  \draw[edge] (Gf)--(Mf) (Gf2)--(Mf2) (Gf12)--(Mf12);
  \draw[hidden] (Gf1)--(Mf1);
  \draw[hidden] (Gf)--(Gf1)--(Gf12);
  \draw[edge] (Mf)--(Mf1) (Gf2)--(Gf12) (Mf2)--(Mf12);
  \draw[edge] (Gf)--(Gf2) (Mf)--(Mf2) (Mf1)--(Mf12);
  \draw[edge] (f)--(Gf) (f2)--(Gf2) (f12)--(Gf12) (f1)--(Gf1);
  \draw[edge] (s)--(Mf) (s1)--(Mf1) (s2)--(Mf2) (s12)--(Mf12);

  \foreach \p in {Gf,Mf,Gf2,Mf2,Gf1,Mf1,Gf12,Mf12,f,s,f2,s2,f1,s1,f12,s12}
    \node[vtx] at (\p) {};

  \node[lab, above left]  at (Gf12) {$G_\mu(f^{12})$};
  \node[lab, above right] at (Mf12) {$\mu(f^{12*})$};
  \node[lab, above left]  at (Gf2)  {$G_\mu(f^{2})$};
  \node[lab, above left]  at (Mf2)  {$\mu(f^{2*})$};
  \node[lab, below right] at (f12)  {$f^{12}$};
  \node[lab, below right] at (s12)  {$f^{12*}$};
  \node[lab, below left]  at (f2)   {$f^{2}$};
  \node[lab, below left]  at (s2)   {$f^{2*}$};
  \node[lab, above right] at (f1)   {$f^{1}$};
  \node[lab, above right] at (s1)   {$f^{1*}$};
  \node[lab, above left]  at (f)    {$f$};
  \node[lab, above left]  at (s)    {$f^{*}$};
  \node[lab, below right] at (Gf1)  {$G_\mu(f^{1})$};
  \node[lab, below right] at (Mf1)  {$\mu(f^{1*})$};
  \node[lab, below left]  at (Gf)   {$G_\mu(f)$};
  \node[lab, below right] at (Mf)   {$\mu(f^{*})$};

  \node at (3.02,-0.95) {(a) The four-dimensional cube.};
\end{scope}

\begin{scope}[shift={(-1.155,-3.31)}]
  \node (a)   at (0,0)        {$f$};
  \node (a1)  at (2.54,0)     {$f^{1}$};
  \node (a2)  at (0,1.7)      {$f^{2}$};
  \node (a12) at (2.54,1.7)   {$f^{12}$};
  \node (b)   at (5.76,0)     {$f^{*}$};
  \node (b1)  at (8.30,0)     {$f^{1*}$};
  \node (b2)  at (5.76,1.7)   {$f^{2*}$};
  \node (b12) at (8.30,1.7)   {$f^{12*}$};
  \node (c)   at (0,-3.47)    {$G_\mu(f)$};
  \node (c1)  at (2.54,-3.47) {$G_\mu(f^{1})$};
  \node (c2)  at (0,-1.78)    {$G_\mu(f^{2})$};
  \node (c12) at (2.54,-1.78) {$G_\mu(f^{12})$};
  \node (d)   at (5.76,-3.47) {$\mu(f^{*})$};
  \node (d1)  at (8.30,-3.47) {$\mu(f^{1*})$};
  \node (d2)  at (5.76,-1.78) {$\mu(f^{2*})$};
  \node (d12) at (8.30,-1.78) {$\mu(f^{12*})$};

  \foreach \v in {a,b,c,d} {
    \draw[sq] (\v)   -- node[above]{$D_{t_1}$} (\v1);
    \draw[sq] (\v2)  -- node[below]{$D_{t_1}$} (\v12);
    \draw[sq] (\v)   -- node[right]{$D_{t_2}$} (\v2);
    \draw[sq] (\v1)  -- node[left] {$D_{t_2}$} (\v12);
  }

  \draw[link] (3.22,0.85)   -- node[above=2pt]{$*$} (5.04,0.85);
  \draw[link] (3.22,-2.63)  -- node[above=2pt]{$*$} (5.04,-2.63);
  \draw[link] (1.27,-0.47)  -- node[right=3pt]{$G_\mu$} (1.27,-1.27);
  \draw[link] (7.03,-0.47)  -- node[right=3pt]{$\mu$} (7.03,-1.27);

  \node at (4.17,-4.22) {(b) The four Darboux squares and their vertexwise maps.};
\end{scope}
\end{tikzpicture}
\caption{The four-dimensional permutability of Theorem~\ref{p1:diagram}. (a) Two nested cubes with corresponding vertices joined. (b) The four Darboux squares. The double arrows act on all four vertices: from left to right by Christoffel duality, and from top to bottom by $G_\mu$ on the left and by $\mu$ on the right.}
\label{p1:permutability}
\end{figure}

The surface $f^{12}$ is thus obtained without solving another Riccati equation. For $\mu=\iota$, the M\"obius square consists of the stereographic images of the four Gauss maps; its Christoffel dual, the Goursat square, supplies the planar data for the Weierstrass representation~\eqref{p1:weierstrass}. Section~\ref{sec:liouville} shows that the conformal factor of $f^{12}$ satisfies the Liouville superposition formula of~\citep[Theorem 3.3]{CFT}.

In stereographic coordinates, the cross ratio in~\eqref{p1:mobiusface} reads
\[
\operatorname{cr}(g,g^1,g^{12},g^2)=\frac{t_2}{t_1}.
\]
Solving for the remaining vertex yields
\begin{equation}
g^{12}=g+(t_2-t_1)\frac{(g^1-g)(g^2-g)}{t_2(g^2-g)-t_1(g^1-g)}.
\label{p1:superposition}
\end{equation}
Indeed,
\[
\frac1{g^{12}-g}=\frac1{t_2-t_1}
\left(\frac{t_2}{g^1-g}-\frac{t_1}{g^2-g}\right).
\]
Differentiating this identity and using $(1/(g^i-g))_z=-t_i/g_z+g_z/(g^i-g)^2$, which follows from the Riccati equations $g_z^i=t_i(g^i-g)^2/g_z$, we obtain
\begin{equation}
g_z^{12}=t_1\frac{(g^{12}-g^2)^2}{\partial_z g^2}
=t_2\frac{(g^{12}-g^1)^2}{\partial_z g^1}
=t_1t_2g_z
\left(\frac{g^1-g^2}{t_2(g^2-g)-t_1(g^1-g)}\right)^2.
\label{p1:superderivative}
\end{equation}
Here $g_z^i=\partial_z g^i$; the identities extend meromorphically through removable chart singularities.

\begin{remark}[Parameter signs]\label{p1:signs}
The computation imposes no restriction on the sign of $t_1t_2$. Since $t_2/t_1$ is real and differs from~1, the formula applies to parameters of either the same or opposite signs.
\end{remark}

\begin{remark}[Equal parameters]\label{p1:equal}
For $t_1=t_2$ and distinct first transforms, direct substitution in~\eqref{p1:superposition} returns $g^{12}=g$ wherever the denominator is nonzero. A nontrivial equal-parameter construction can instead be obtained by a limiting procedure, following the idea of Cho, Leschke and Ogata~\citep{CLO}. We leave the corresponding minimal case for future work and assume throughout that $t_1\ne t_2$.
\end{remark}

\subsection{Ends of the superposed minimal surface}\label{sec:ends}
We now locate the branch points of $g^{12}$; by \eqref{p1:fundamental}, these are the points where the metric of $f^{12}$ blows up. To include poles, use the spherical derivative
\[
\gamma^\sharp:=\frac{2|\gamma_z|}{1+|\gamma|^2}.
\]
It extends continuously across poles and vanishes exactly at the branch points of $\gamma:\Sigma\to\mathbb{CP}^1$; its vanishing order is the order of the derivative in a finite spherical coordinate and equals $k-1$ at a pole of order $k$. For the corresponding normalized metric, $e^\omega=1/\gamma^\sharp$, so the metric is regular exactly where $\gamma^\sharp\ne0$. Since rotations of $\mathbb R^3$ preserve the Hopf differential and Darboux pairs~\citep{HJ} and act on the Gauss maps by rotations of $S^2\cong\mathbb{CP}^1$, all statements below are invariant under simultaneous rotations; we use this to make prescribed values of the Gauss maps at a point finite, or one of them zero.

\begin{lemma}\label{p1:branchlemma}
Let $(f,\hat f)$ be a minimal Darboux pair with Gauss maps $g,\hat g$ as in \eqref{p1:minriccati}--\eqref{p1:relative}, and let $g$ be unbranched at $z_0$. Then $\hat g$ is branched at $z_0$ if and only if $\hat g(z_0)=g(z_0)$ in $\mathbb{CP}^1$. In particular, every pole of $\hat g$ at which $g$ is finite and unbranched is simple.
\end{lemma}
\begin{proof}
By Proposition~\ref{p1:cdprop} (see \eqref{p1:relative}),
\[
\hat f-f=\frac1t(\hat\nu-\nu)^{-1},
\]
and, since $\nu=\iota(g\mathbf j)$ and $\hat\nu=\iota(\hat g\mathbf j)$, we have $\hat g(z_0)=g(z_0)$ in $\mathbb{CP}^1$ if and only if $\hat\nu(z_0)=\nu(z_0)$.

If $\hat\nu(z_0)\ne\nu(z_0)$, then $\hat f-f$ is finite and nonzero near $z_0$. Since $df^*=d\nu$ has rank two, \eqref{p1:qriccati},
\[
d\hat f=t(\hat f-f)\,d\nu\,(\hat f-f),
\]
shows that $d\hat f$ also has rank two, so $\hat f$ is immersed at $z_0$. Hence its conformal factor $e^{\hat\omega}=1/\hat g^\sharp$ is finite at $z_0$, that is, $\hat g^\sharp(z_0)\ne0$.

If $\hat\nu(z_0)=\nu(z_0)$, then $|\hat f-f|\to\infty$ as $z\to z_0$. If $\hat g^\sharp(z_0)\ne0$, the conformal factor $1/\hat g^\sharp$ would be bounded near $z_0$, and so would $\hat f$; hence $\hat g^\sharp(z_0)=0$.

Finally, at a pole of $\hat g$ where $g$ is finite, we have $\hat g\ne g$, so $\hat g$ is unbranched there, and the pole is simple.
\end{proof}

\begin{corollary}\label{p1:regularcor}
If $g$ is unbranched, then $g^{12}$ is unbranched wherever $g^1\ne g^2$ in $\mathbb{CP}^1$. Every pole at such a point is simple, and the metric of $f^{12}$ is regular there.
\end{corollary}
\begin{proof}
Work in a spherical chart in which $g,g^1,g^2$ are finite. Set $\Delta:=t_2(g^2-g)-t_1(g^1-g)$. Equation~\eqref{p1:superposition} gives
\begin{equation}
g^{12}-g^1=\frac{t_1(g^1-g)(g^1-g^2)}{\Delta},\qquad
g^{12}-g^2=\frac{t_2(g^2-g)(g^1-g^2)}{\Delta}.
\label{p1:superdifferences}
\end{equation}
If $g^1\ne g$, then $g^1$ is unbranched by Lemma~\ref{p1:branchlemma}, and $g^{12}\ne g^1$, interpreted projectively if $\Delta=0$. Apply Lemma~\ref{p1:branchlemma} to $(g^1,g^{12})$. Otherwise $g^2\ne g^1=g$, and the same argument applies to $(g^2,g^{12})$.
\end{proof}

\begin{proposition}\label{p1:branchorders}
Let $g$ be unbranched and $t_1\ne t_2$. Then $g^{12}$ is branched exactly where $g^1=g^2$ in $\mathbb{CP}^1$. At such a point $z_0$:
\begin{enumerate}
\item[(a)] if $g^1(z_0)=g^2(z_0)\ne g(z_0)$, then $(g^{12})^\sharp$ has a zero of order exactly $2$;
\item[(b)] if $g^1(z_0)=g^2(z_0)=g(z_0)$, then $(g^{12})^\sharp$ has a zero of order exactly $4$.
\end{enumerate}
\end{proposition}
\begin{proof}
Corollary~\ref{p1:regularcor} treats unequal first values. Rotate so that the three initial values at $z_0$ are finite; then $g_z(z_0)\ne0$, since $g$ is unbranched. In (a), $\Delta(z_0)=(t_2-t_1)(g^1(z_0)-g(z_0))\ne0$, and
\[
(g^1-g^2)_z(z_0)=\frac{(t_1-t_2)(g^1(z_0)-g(z_0))^2}{g_z(z_0)}\ne0.
\]
Thus \eqref{p1:superderivative} shows that the derivative vanishes to order $2$, and \eqref{p1:superposition} gives $g^{12}(z_0)=g^1(z_0)$. In (b), put $p_0=g_z(z_0)$. The Riccati equations give
\[
g^i-g=-\int_{z_0}^z g_z\,dz+\frac{t_ip_0}{3}(z-z_0)^3+O((z-z_0)^4),
\]
so
\[
\begin{aligned}
g^1-g^2&=\frac{(t_1-t_2)p_0}{3}(z-z_0)^3+O((z-z_0)^4),\\
\Delta&=-(t_2-t_1)p_0(z-z_0)+O((z-z_0)^2).
\end{aligned}
\]
Therefore \eqref{p1:superderivative} shows that the derivative vanishes to order exactly $4$, and \eqref{p1:superposition} gives $g^{12}(z_0)=g(z_0)$.
\end{proof}

Proposition~\ref{p1:branchorders} requires no genericity assumption on the integration constants. We next determine the multiplicities of the corresponding ends and compare the result with~\citep{MRT} in Remark~\ref{p1:comparison}.

\begin{proposition}[cf. \citep{JM,MRT}]\label{p1:planarprop}
Let $f$ have Weierstrass data $(\gamma,-dz/\gamma_z)$ and be defined and single-valued on a punctured neighborhood of $z_0$, with $\gamma$ meromorphic near $z_0$. If $\gamma^\sharp$ has a zero of order $m\ge1$, then $m\ge2$, and $z_0$ is a complete planar end of multiplicity $m-1$. It is embedded if and only if $m=2$.
\end{proposition}
\begin{proof}
After a rotation of $\mathbb R^3$ (see the paragraph before Lemma~\ref{p1:branchlemma}), we may assume that $\gamma(z_0)=0$. Then, with $c_0\ne0$,
\[
\gamma=c_0(z-z_0)^{m+1}\bigl(1+O(z-z_0)\bigr),\qquad
h=-\frac{(z-z_0)^{-m}}{(m+1)c_0}\bigl(1+O(z-z_0)\bigr),
\]
so that $\gamma h\,dz$ and $\gamma^2h\,dz$ are holomorphic at $z_0$. By \eqref{p1:weierstrass},
\[
f_1+\mathrm if_2=\frac12\left(\overline{\int^z h\,dz}-\int^z\gamma^2h\,dz\right),\qquad
f_3=\operatorname{Re}\int^z\gamma h\,dz.
\]
Since $f$ is single-valued, integration along a small circle around $z_0$ gives
\[
0=\oint d(f_1+\mathrm if_2)=\frac12\overline{\oint h\,dz}
=-\pi\mathrm i\,\overline{\operatorname{Res}_{z_0}(h\,dz)}.
\]
For $m=1$ this residue would be $-1/(2c_0)\ne0$; hence $m\ge2$. Moreover, all three coordinate forms in \eqref{p1:weierstrass} have zero residue at $z_0$. Integrating term by term, we obtain
\[
f_1+\mathrm if_2=\overline{\frac{(z-z_0)^{1-m}}{2(m^2-1)c_0}\bigl(1+O(z-z_0)\bigr)},\qquad
f_3=O(1),
\]
while \eqref{p1:fundamental} gives
\[
ds=\frac{|dz|}{\gamma^\sharp}\asymp|z-z_0|^{-m}|dz|.
\]
Thus $z_0$ is a complete planar end of multiplicity $m-1$, which is embedded exactly when $m=2$ by~\citep{JM}.
\end{proof}

\begin{corollary}\label{p1:endcor}
If $g$ is unbranched at $z_0$ and $g^1(z_0)=g^2(z_0)$, then $f^{12}$ is single-valued near the puncture and has a complete planar end there: an embedded end in case (a), and an end of multiplicity $3$ in case (b). All componentwise residues of its Weierstrass integrand vanish.
\end{corollary}
\begin{proof}
Proposition~\ref{p1:cdprop} gives, on a sufficiently small punctured disk,
\[
f^{12}=f+\frac1{t_1}(f^{1*}-f^*)^{-1}+\frac1{t_2}(f^{12*}-f^{1*})^{-1}.
\]
The right-hand side is single-valued. Apply Propositions~\ref{p1:branchorders} and~\ref{p1:planarprop} with $m=2$ or $4$. After rotating the final Gauss value at $z_0$ to zero, we see from the preceding proof that all three residues vanish. The inverse rotation preserves this vanishing.
\end{proof}

\begin{remark}[Comparison with \citep{MRT}]\label{p1:comparison}
The minimal surfaces in this subsection have Hopf differential $dz^2/2$ on their regular domains and hence have no umbilics there. In case~(a), $f^1$ is regular at $z_0$, so~\citep[Theorem 3.5]{MRT} already implies that the new end is embedded and planar. In case~(b), both intermediate surfaces have embedded planar ends at $z_0$. Theorem~3.8 of~\citep{MRT} allows the transformed end to be planar or removable; our computation gives a multiplicity of exactly~3. Thus superposition produces a non-embedded planar end over a regular umbilic-free point of the seed. This differs from the mechanism involving seed umbilics in~\citep[Theorem 3.6]{MRT}. Section~\ref{sec:examples} gives explicit examples.
\end{remark}

\subsection{Recovering the superposition principle for the Liouville equation}\label{sec:liouville}
For a minimal surface $f$ with Hopf differential $dz^2/2$ and metric $I_f=e^{2\omega}(du^2+dv^2)$, the Gaussian curvature is $K=-e^{-4\omega}$. The Gauss equation therefore reduces to the Liouville equation
\[
\omega_{uu}+\omega_{vv}=e^{-2\omega}.
\]
By \eqref{p1:fundamental}, the conformal factors of $f$, $f^i$ and $f^{12}$ are
\[
e^\omega=\frac{1+|g|^2}{2|g_z|},\qquad
e^{\omega^i}=\frac{1+|g^i|^2}{2|g_z^i|}\quad(i=1,2),\qquad
e^{\omega^{12}}=\frac{1+|g^{12}|^2}{2|g_z^{12}|}.
\]
The metric functions $\omega$, $\omega^i$ and $\omega^{12}$ solve the Liouville equation on their respective regular domains. Substituting~\eqref{p1:superposition}--\eqref{p1:superderivative} into the expression for $e^{\omega^{12}}$ gives the following superposition principle for the Liouville equation.
\begin{corollary}\label{p1:liouville}
Let $g^i=D_{t_i}(g)$ $(i=1,2)$ be Darboux transforms with distinct nonzero real parameters $t_1,t_2$. Let $\omega$ and $\omega^i$ denote the normalized metric functions of $f$ and $f^i$, respectively. Then
\[
e^{\omega^{12}}=e^\omega\,
\frac{|t_1|e^{\omega^2}+|t_2|e^{\omega^1}-\operatorname{sgn}(t_1t_2)B}
{|t_2|e^{\omega^2}+|t_1|e^{\omega^1}-B},\qquad
B=|g_z|\operatorname{Re}\frac{1+g^1\overline{g^2}}{(g^1-g)\overline{(g^2-g)}}.
\]
Where the conformal factor is finite and positive, its logarithm $\omega^{12}$ solves the Liouville equation. The formula is initially evaluated where all four metric functions are regular and the displayed denominators are nonzero, and extends by continuation wherever $f^{12}$ is regular.
\end{corollary}
\begin{proof}
The metric functions are invariant under a simultaneous rotation of the surfaces. So is $B$, since
\[
B=2e^{-\omega}\frac{\langle\nu^1-\nu,\nu^2-\nu\rangle}{|\nu^1-\nu|^2|\nu^2-\nu|^2}.
\]
We may therefore rotate the surfaces so that $g=0$ at the point under consideration. Set $\beta_i=g^i$ $(i=1,2)$, $p=|g_z|$, $\delta=t_2-t_1$ and $\Delta=t_2\beta_2-t_1\beta_1$ as in the proof of Corollary~\ref{p1:regularcor}. The Riccati equations give
\[
e^{\omega^1}=\frac{p(1+|\beta_1|^2)}{2|t_1||\beta_1|^2},\qquad
e^{\omega^2}=\frac{p(1+|\beta_2|^2)}{2|t_2||\beta_2|^2},\qquad
B=p\operatorname{Re}\frac{1+\beta_1\overline{\beta_2}}{\beta_1\overline{\beta_2}}.
\]
Consequently the denominator and numerator of the quotient in the statement simplify to
\[
\begin{aligned}
|t_2|e^{\omega^2}+|t_1|e^{\omega^1}-B
&=\frac{p|\beta_1-\beta_2|^2}{2|\beta_1\beta_2|^2},\\
|t_1|e^{\omega^2}+|t_2|e^{\omega^1}-\operatorname{sgn}(t_1t_2)B
&=\frac{p\bigl(|\Delta|^2+\delta^2|\beta_1\beta_2|^2\bigr)}{2|t_1t_2||\beta_1\beta_2|^2}.
\end{aligned}
\]
Since $e^\omega=1/(2p)$, their quotient yields
\[
e^{\omega^{12}}=\frac{|\Delta|^2+\delta^2|\beta_1\beta_2|^2}{2|t_1t_2|p|\beta_1-\beta_2|^2}
=\frac{1+|g^{12}|^2}{2|g_z^{12}|}
\]
by \eqref{p1:superposition}--\eqref{p1:superderivative}. Theorem~\ref{p1:minimal} and the Gauss equation prove the assertion.
\end{proof}

To compare Corollary~\ref{p1:liouville} with~\citep[Theorem 3.3]{CFT}, recall that the latter assumes two special Ribaucour transforms of a common seed, equivalently minimal Darboux transforms. It does not apply to arbitrary triples of Liouville solutions. On a common regular patch, take $\phi=\omega$, the principal frame $e_1=e^{-\omega}f_u$, $e_2=e^{-\omega}f_v$, and the auxiliary functions $\Omega_i,\Omega_{i,1},\Omega_{i,2},W_i$ of~\citep[(41)]{CFT}. Then
\[
(W_i)_u=e^{-\omega}\Omega_{i,1},\qquad
(W_i)_v=-e^{-\omega}\Omega_{i,2},\qquad
\Omega_{i,1}^2+\Omega_{i,2}^2+W_i^2=2c_i\Omega_iW_i.
\]
With $X_i=\nu^i-\nu$ and $L_i=|X_i|^2$, Equations (12)--(13) of~\citep{CFT} give
\[
X_i=\frac{\Omega_{i,1}e_1+\Omega_{i,2}e_2-W_i\nu}{c_i\Omega_i},\qquad
f^i-f=-\frac{\Omega_i}{W_i}X_i=-\frac2{c_i}\frac{X_i}{L_i},
\]
and comparison with \eqref{p1:relative} yields $c_i=2t_i$. Using the derivatives of $W_i$, the expression for $X_i$ and $\langle X_i,\nu\rangle=-L_i/2$, we can write the auxiliary quantity in~\citep{CFT}, denoted here by $\mathcal B$, as
\[
\begin{aligned}
\mathcal B
&=e^{-\omega}+e^\omega\langle\nabla\log|W_1|,\nabla\log|W_2|\rangle\\
&=e^{-\omega}\left(1+\frac{\Omega_{1,1}\Omega_{2,1}+\Omega_{1,2}\Omega_{2,2}}{W_1W_2}\right)
=\frac{4e^{-\omega}\langle X_1,X_2\rangle}{L_1L_2}=2B,
\end{aligned}
\]
where $\nabla=(\partial_u,\partial_v)$; the last equality is the expression of $B$ in the proof of Corollary~\ref{p1:liouville}. With positive metric factors $e^{\omega^i}=e^{-\omega}|\Omega_i/W_i|$, the formula of~\citep[Theorem 3.3]{CFT} reads
\[
e^{\omega^{12}}=e^\omega\frac{|c_1|e^{\omega^2}+|c_2|e^{\omega^1}-\operatorname{sgn}(c_1c_2)\mathcal B}
{|c_1|e^{\omega^1}+|c_2|e^{\omega^2}-\mathcal B},
\]
which is exactly Corollary~\ref{p1:liouville} under $c_i=2t_i$ and $\mathcal B=2B$. For $c_1,c_2>0$ this gives another local proof of~\citep[Theorem 3.3]{CFT}. For $c_i<0$ the quadratic constraint gives $\operatorname{sgn}(\Omega_i/W_i)=\operatorname{sgn}c_i$, so the conformal factor in~\citep[Theorem 3.2]{CFT} has to be read with absolute values. In particular, the minimal surface determined by the solution in~\citep[Theorem 3.3]{CFT} is, up to congruence, the common Darboux transform $f^{12}$. All comparisons are made away from the zeros of the auxiliary functions and denominators, with continuation on regular parts.

\section{Examples}\label{sec:examples}
We apply~\eqref{p1:superposition}--\eqref{p1:superderivative} to the explicit minimal Darboux transforms of the Enneper surface and the catenoid. For the catenoid, a suitable parameter choice gives transforms on $\mathbb C^*$, yielding complete minimal surfaces of finite total curvature (Proposition~\ref{p1:global}). Furthermore, the subfamily satisfying~\eqref{p1:cat-normalization} realizes a L\'opez--Ros deformation (Proposition~\ref{p1:lopezros}).

We begin with the Enneper surface and determine when its superposed transform has a planar end of multiplicity~3.

\begin{example}[Enneper surface]\label{p1:enneper}

In our normalization, the Enneper surface has Gauss map
\[
g(z)=z\qquad(z\in\mathbb C).
\]
A Darboux transform $\hat g$ of $g$ is a solution to the ordinary differential equation
\[
\hat g_z=t(\hat g-z)^2.
\]
Its solutions are

\[
\hat g=z+\frac{1+Ce^{2\sqrt t\,z}}{\sqrt t\,(1-Ce^{2\sqrt t\,z})}
\qquad(C\in\mathbb C\cup\{\infty\}).
\]
In the rest of this example we assume that $t>0$. Then
\[
\hat g_z=t(\hat g-z)^2=
\frac{(1+Ce^{2\sqrt t\,z})^2}{(1-Ce^{2\sqrt t\,z})^2}.
\]
For $C\in\mathbb C^*$, its zeros satisfy $Ce^{2\sqrt t\,z}=-1$.
If $C=|C|e^{\mathrm i\theta}$, these zeros are
\[
z=\frac{\log(1/|C|)}{2\sqrt t}
+\mathrm i\frac{(\pi-\theta)+2k\pi}{2\sqrt t},\qquad k\in\mathbb Z.
\]
They are zeros of order~2, arranged along a line parallel to the imaginary axis
with spacing $\pi/\sqrt t$. The corresponding points of $\hat f$ are embedded
planar ends (Proposition~\ref{p1:planarprop}; cf.~\citep[Theorem~3.5]{MRT}).
Single-valuedness follows from Proposition~\ref{p1:cdprop}.

Write $g^i=D_{t_i}(g)$ $(i=1,2)$ as
\[
g^i=z+\frac{2(1+C_{\lambda_i}e^{\lambda_i z})}
{\lambda_i(1-C_{\lambda_i}e^{\lambda_i z})},
\qquad \lambda_i:=2\sqrt{t_i},\quad t_1\ne t_2.
\]
By \eqref{p1:superposition}, $g^{12}$ is given algebraically by
\[
g^{12}=z+\frac{2(\lambda_2^2-\lambda_1^2)}{\lambda_1\lambda_2}
\frac{\mathcal P(z)}{\mathcal D(z)},
\]
where
\begin{align*}
\mathcal P(z)&:=C_{\lambda_1}C_{\lambda_2}e^{(\lambda_1+\lambda_2)z}
+C_{\lambda_1}e^{\lambda_1z}+C_{\lambda_2}e^{\lambda_2z}+1,\\
\mathcal D(z)&:=(\lambda_1-\lambda_2)C_{\lambda_1}C_{\lambda_2}e^{(\lambda_1+\lambda_2)z}
-(\lambda_1+\lambda_2)C_{\lambda_1}e^{\lambda_1z}+(\lambda_1+\lambda_2)C_{\lambda_2}e^{\lambda_2z}-(\lambda_1-\lambda_2).
\end{align*}
Using \eqref{p1:superderivative}, we obtain
\[
g_z^{12}=\left(\frac{\mathcal R(z)}{\mathcal D(z)}\right)^2,
\]
where
\begin{align*}
\mathcal R(z):=(\lambda_1-\lambda_2)C_{\lambda_1}C_{\lambda_2}e^{(\lambda_1+\lambda_2)z}
+(\lambda_1+\lambda_2)C_{\lambda_1}e^{\lambda_1z}-(\lambda_1+\lambda_2)C_{\lambda_2}e^{\lambda_2z}-(\lambda_1-\lambda_2).
\end{align*}
Since $g=z$ is unbranched, Proposition~\ref{p1:branchorders} and Corollary~\ref{p1:endcor} identify all ends of the single-valued surface $f^{12}$ in $\mathbb C$: they occur where $g^1=g^2$ and are embedded planar ends unless $g^1=g^2=z$, when their multiplicity is~3. The latter points are the common zeros of $g_z^1$ and $g_z^2$. They occur precisely when $C_{\lambda_i}=-e^{-\lambda_i z_0}$ for $i=1,2$ and some $z_0\in\mathbb C$. If $\lambda_1/\lambda_2=p/q$ with coprime positive integers $p,q$, these points are $z_0+2\pi\mathrm i np/\lambda_1$, $n\in\mathbb Z$; if the ratio is irrational, $z_0$ is the only such point (see Figure~\ref{p3-Enneper}).

\begin{figure}[H]
\centering
\includegraphics[width=0.8\linewidth]{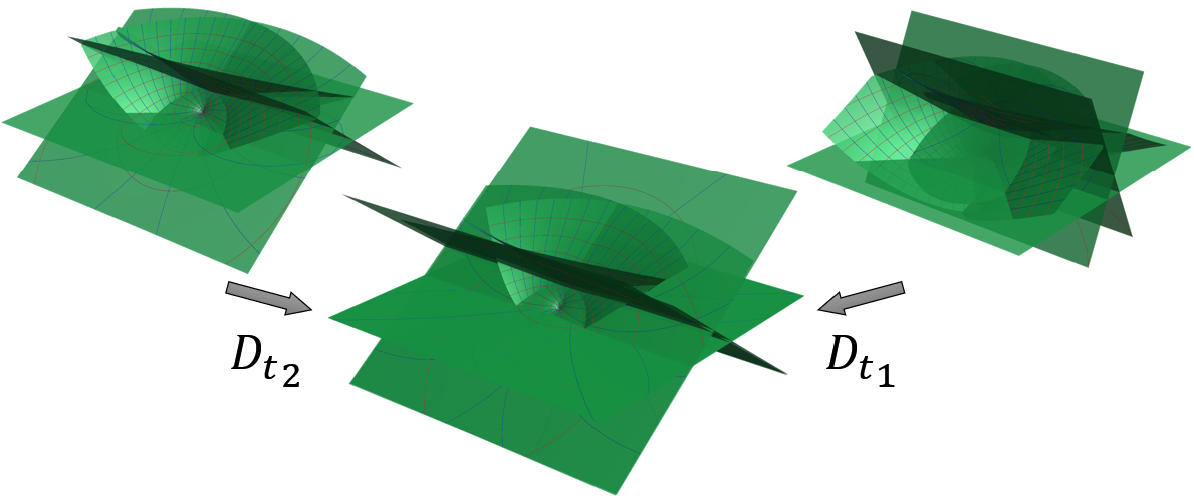}
\caption{Superposition principle for minimal Darboux transforms of the Enneper surface. Left and right: the single transforms $f^1=D_{t_1}(f)$ and $f^2=D_{t_2}(f)$. Middle: the superposed surface $f^{12}$. The initial values are chosen as $C_{\lambda_i}=-e^{-\lambda_i z_0}$ $(i=1,2)$ for some $z_0\in\mathbb C$. At $z_0$, both $f^1$ and $f^2$ have embedded planar ends, whereas $f^{12}$ has a planar end of multiplicity~3. The parameters are $\lambda_1 = 2$ and $\lambda_2=4$, with $C_{\lambda_i}=-1 \ (i=1,2)$.}\label{p3-Enneper}
\end{figure}
\end{example}

For the catenoid, we first locate the planar ends of a single periodic Darboux transform.

\begin{example}[Catenoid]\label{p1:catenoid}

The catenoid is given by the Gauss map
\[
g(z)=e^z\qquad(z\in\mathbb C),
\]
i.e. by the classical Weierstrass data
\[
(g,h\,dz)=(e^z,-e^{-z}dz)=\left(\zeta,-\frac{d\zeta}{\zeta^2}\right),
\qquad \zeta:=e^z.
\]
For a complete minimal surface of finite total curvature, we use
Osserman's Gauss-map formula \citep{Osserman} and the Jorge--Meeks
end formula \citep{JM}:
\[
\int_M K\,dA=-4\pi\deg g
=2\pi\left(\chi(\overline M)-\sum_j(I_j+1)\right),
\]
where $\overline M$ is the compactification and $I_j$ is the end
multiplicity; an end is embedded exactly when $I_j=1$.

A Darboux transform $\hat g$ satisfies
\[
\hat g_z=t\frac{(\hat g-e^z)^2}{e^z}.
\]
For $t>-1/4$, $t\ne0$, set $\lambda=\sqrt{1+4t}$. The solutions have the following form, where $C_\lambda\in\mathbb C\cup\{\infty\}$ and the value at $C_\lambda=\infty$ is understood as a limit:
\begin{equation}
\hat g(z)=e^z\left\{1+\frac{2}{\lambda^2-1}
\left(1+\lambda\frac{1+C_\lambda e^{\lambda z}}
{1-C_\lambda e^{\lambda z}}\right)\right\}.\label{p1:cat-single}
\end{equation}
From now on, we choose $t$ so that $\lambda\in\mathbb Z_{\ge2}$. Then $\hat g$ is
$2\pi\mathrm i$-periodic and defines a meromorphic function of $\zeta=e^z$ on $\mathbb C^*$.
For $C_\lambda\in\mathbb C^*$, the zeros of $\hat g_z$ satisfy
\[
C_\lambda e^{\lambda z}=-\frac{\lambda+1}{\lambda-1}.
\]
With $\theta=\arg C_\lambda$, these zeros are
\begin{equation}
z=\frac1\lambda\log\left(\frac{\lambda+1}{\lambda-1}|C_\lambda|^{-1}\right)
+\mathrm i\frac{\pi-\theta+2k\pi}{\lambda},\qquad k\in\mathbb Z.
\label{p1:cat-singlezeros}
\end{equation}
There are $\lambda$ such points modulo $2\pi\mathrm i$. For integer $\lambda\geq2$, the corresponding surface is one of the surfaces $\tilde X_{(\lambda,1)}$ in~\citep[Corollary~2.9]{CFT}; see also~\citep[Section~4.1]{MRT}. It is complete on $\mathbb C^*$ minus these $\lambda$ embedded planar ends, with catenoidal ends at $0,\infty$ and total curvature $-4\pi(\lambda+1)$. More generally, the corresponding parameter $\sqrt{1-2c}$ in~\citep[Propositions~2.4(c) and~2.8]{CFT} may be any positive rational number $n/m\ne1$ in lowest terms. The two nonplanar ends then have geometric index~$m$; here $m=1$.
\end{example}

The following remark relates the planar ends of a single transform to the lobes of a bubbleton.

\begin{remark}\label{p1:bubbletons}
For a surface of revolution, the resonance parameters in~\citep{CLO}, at which
all Darboux transforms close up, are $\varrho_k=(k^2-1)/4$ with
$k\in\mathbb Z_{>1}$. Closed transforms that are not surfaces of revolution
are called bubbletons and have $k$ lobes. The Christoffel and Riccati
normalizations in~\citep{CLO} agree with ours, so $t=(\lambda^2-1)/4$,
$\lambda\in\mathbb Z_{\ge2}$, is a resonance parameter with $k=\lambda$.
Thus the minimal catenoid transforms with $C_\lambda\in\mathbb C^*$ are
bubbletons in this sense, and their $\lambda$ planar ends in
\eqref{p1:cat-singlezeros} correspond to the lobes.
\end{remark}

To superpose two such transforms, take $\lambda_1\ne\lambda_2$ and $C_{\lambda_1},C_{\lambda_2}\in\mathbb C^*$.
The choices $C_\lambda=0,\infty$ give constant multiples of $g$.
Write $g^i=D_{t_i}(g)$ as
\begin{equation}
g^i(z)=e^z\left\{1+\frac{2}{\lambda_i^2-1}
\left(1+\lambda_i\frac{1+C_{\lambda_i}e^{\lambda_i z}}
{1-C_{\lambda_i}e^{\lambda_i z}}\right)\right\},
\quad\lambda_i=\sqrt{1+4t_i}.
\label{p1:cat-pair}
\end{equation}
The zeros of $g_z^i$ are
\begin{equation}
z^i=\frac1{\lambda_i}\log\left(\frac{\lambda_i+1}{\lambda_i-1}
|C_{\lambda_i}|^{-1}\right)
+\mathrm i\frac{\pi-\arg C_{\lambda_i}+2k_i\pi}{\lambda_i},
\quad k_i\in\mathbb Z.
\label{p1:cat-pairzeros}
\end{equation}

The function $g^{12}$ determined by \eqref{p1:superposition} is
\begin{equation}
g^{12}=\frac{\zeta P^1(\zeta)}{(\lambda_1^2-1)(\lambda_2^2-1)D(\zeta)},\label{p1:cat-super}
\end{equation}
where
\begin{align*}
P^1(\zeta)&:=(\lambda_1-1)^2(\lambda_2-1)^2(\lambda_1-\lambda_2)
C_{\lambda_1}C_{\lambda_2}\zeta^{\lambda_1+\lambda_2}\\
&\quad-(\lambda_1-1)^2(\lambda_2+1)^2(\lambda_1+\lambda_2)
C_{\lambda_1}\zeta^{\lambda_1}\\
&\quad+(\lambda_1+1)^2(\lambda_2-1)^2(\lambda_1+\lambda_2)
C_{\lambda_2}\zeta^{\lambda_2}
-(\lambda_1+1)^2(\lambda_2+1)^2(\lambda_1-\lambda_2),\\
D(\zeta)&:=(\lambda_1-\lambda_2)C_{\lambda_1}C_{\lambda_2}\zeta^{\lambda_1+\lambda_2}
-(\lambda_1+\lambda_2)C_{\lambda_1}\zeta^{\lambda_1}+(\lambda_1+\lambda_2)C_{\lambda_2}\zeta^{\lambda_2}-(\lambda_1-\lambda_2).
\end{align*}
Furthermore, by \eqref{p1:superderivative},
\begin{equation}
g_z^{12}=\frac{\zeta}{(\lambda_1^2-1)(\lambda_2^2-1)}
\left(\frac{P^2(\zeta)}{D(\zeta)}\right)^2,\label{p1:cat-superderivative}
\end{equation}
where
\begin{align*}
P^2(\zeta)&:=(\lambda_1-1)(\lambda_2-1)(\lambda_1-\lambda_2)
C_{\lambda_1}C_{\lambda_2}\zeta^{\lambda_1+\lambda_2}\\
&\quad+(\lambda_1-1)(\lambda_2+1)(\lambda_1+\lambda_2)
C_{\lambda_1}\zeta^{\lambda_1}\\
&\quad-(\lambda_1+1)(\lambda_2-1)(\lambda_1+\lambda_2)
C_{\lambda_2}\zeta^{\lambda_2}
-(\lambda_1+1)(\lambda_2+1)(\lambda_1-\lambda_2).
\end{align*}
Since $P^1(0)\ne0\ne D(0)$ and $\deg P^1=\deg D=\lambda_1+\lambda_2$, as a
meromorphic function of $\zeta$ we have
\[
\deg(g^{12})=\lambda_1+\lambda_2+1-\ell,
\]
where $\ell$ is the degree of the greatest common divisor of $P^1(\zeta)$ and $D(\zeta)$.

\subsection{Global properties of superposed catenoids}\label{sec:catenoid-global}
To compute the degree of $g^{12}$, we first determine when $P^1$ and $D$ have a common root.

\begin{lemma}\label{p1:commonroots}
$P^1(\alpha)=D(\alpha)=0$ for some $\alpha\in\mathbb C^*$ if and only if
$g^1(z_0)=g^2(z_0)=g(z_0)$ with $e^{z_0}=\alpha$, that is, if and only if
\begin{equation}
C_{\lambda_i}=-\frac{\lambda_i+1}{\lambda_i-1}\alpha^{-\lambda_i}
\qquad(i=1,2).\label{p1:cat-normalization}
\end{equation}
Furthermore, every common root of $P^1$ and $D$ is a simple root of $D$.
\end{lemma}

\begin{proof}
Set $y_i=C_{\lambda_i}\zeta^{\lambda_i}$ and write $g^i=g(1+s_i)$. By \eqref{p1:superposition},
\[
\frac{g^{12}}g=\frac{t_2s_2-t_1s_1+(t_2-t_1)s_1s_2}{t_2s_2-t_1s_1}.
\]
Multiplying the numerator and denominator by $(1-y_1)(1-y_2)$ gives $P^1$ and $D$, respectively, up to nonzero constant factors.
No common root arises at the zeros of these factors. In fact,
\[
D\big|_{y_1=1}=2\lambda_1(y_2-1),\qquad
D\big|_{y_2=1}=-2\lambda_2(y_1-1),
\]
so $D\ne0$ if exactly one of $y_1,y_2$ equals $1$. If $y_1=y_2=1$, one has
$P^1=-8\lambda_1\lambda_2(\lambda_1-\lambda_2)(\lambda_1+\lambda_2)\ne0$.
Away from these zeros, $P^1=D=0$ is therefore equivalent to
\[
t_2s_2=t_1s_1\quad\text{and}\quad s_1s_2=0,
\]
which, since $t_1t_2\ne0$, is equivalent to $s_1=s_2=0$, or $g^1=g^2=g$. By \eqref{p1:cat-single}, this is precisely condition~\eqref{p1:cat-normalization}.
At such a point, $(s_i)_z=t_i s_i^2-s_i-1=-1$, so
\[
(t_2s_2-t_1s_1)_z=t_1-t_2\ne0.
\]
Thus the denominator vanishes to first order.\end{proof}

Once one point with $g^1=g^2=g$ is known, the next proposition counts all such points on $\mathbb C^*$.

\begin{proposition}\label{p1:commoncount}
Assume that~\eqref{p1:cat-normalization} holds. Then the zeros of $g_z^1$ and $g_z^2$ coincide at exactly
$d:=\gcd(\lambda_1,\lambda_2)$ points of $\mathbb C^*$ (i.e. modulo
$2\pi\mathrm i$ in $z$), namely at the solutions of $\zeta^d=\alpha^d$.
\end{proposition}

\begin{proof}
By \eqref{p1:cat-single} and \eqref{p1:cat-normalization}, the zeros of $g_z^i$ satisfy
$(\zeta/\alpha)^{\lambda_i}=1$. By B\'ezout's identity for integers, these two
conditions are equivalent to $(\zeta/\alpha)^d=1$, where
$d=\gcd(\lambda_1,\lambda_2)$. Thus there are exactly $d$ common zeros.
\end{proof}

A change of coordinate places these common zeros at roots of unity and simplifies the polynomials in~\eqref{p1:cat-superderivative}.

\begin{remark}\label{p1:normalizedpolys}
Under condition~\eqref{p1:cat-normalization}, put $X=\zeta/\alpha$. Then
\begin{align*}
\tilde P^1(X)&:=\frac{P^1(\alpha X)}{(\lambda_1+1)(\lambda_2+1)}\\
&=(\lambda_1-1)(\lambda_2-1)(\lambda_1-\lambda_2)X^{\lambda_1+\lambda_2}+(\lambda_1-1)(\lambda_2+1)(\lambda_1+\lambda_2)X^{\lambda_1}\\
&\quad-(\lambda_1+1)(\lambda_2-1)(\lambda_1+\lambda_2)X^{\lambda_2}
-(\lambda_1+1)(\lambda_2+1)(\lambda_1-\lambda_2),\\
\tilde D(X)&:=(\lambda_1-1)(\lambda_2-1)D(\alpha X)\\
&=(\lambda_1+1)(\lambda_2+1)(\lambda_1-\lambda_2)X^{\lambda_1+\lambda_2}+(\lambda_1+1)(\lambda_2-1)(\lambda_1+\lambda_2)X^{\lambda_1}\\
&\quad-(\lambda_1-1)(\lambda_2+1)(\lambda_1+\lambda_2)X^{\lambda_2}
-(\lambda_1-1)(\lambda_2-1)(\lambda_1-\lambda_2),\\
\tilde P^2(X)&:=\frac{P^2(\alpha X)}{(\lambda_1+1)(\lambda_2+1)}\\
&=(\lambda_1-\lambda_2)X^{\lambda_1+\lambda_2}
-(\lambda_1+\lambda_2)X^{\lambda_1}+(\lambda_1+\lambda_2)X^{\lambda_2}-(\lambda_1-\lambda_2),
\end{align*}
and
\[
g^{12}=\frac{\zeta\tilde P^1(X)}{\tilde D(X)}.
\]
These polynomials depend on $\alpha$ only through $X=\zeta/\alpha$.
\end{remark}

By Lemma~\ref{p1:commonroots} and Proposition~\ref{p1:commoncount}, the common roots of $\tilde P^1$ and $\tilde D$ are exactly the $d$th roots of unity. Each is a simple root of $\tilde D$, so their monic greatest common divisor is $X^d-1$. Hence $\ell=d$ and
\[
\deg(g^{12})=\lambda_1+\lambda_2+1-d.
\]
For a meromorphic function $F$, let $\operatorname{ord}_{X_0}F$ denote its order
at $X_0$, positive at a zero and negative at a pole. Since $X=e^z/\alpha$
is a local coordinate, the order is the same in $z$ and $X$.
At each point $X_0^d=1$, Proposition~\ref{p1:branchorders}(b) gives
$\operatorname{ord}_{X_0}g_z^{12}=4$, while Lemma~\ref{p1:commonroots} gives
$\operatorname{ord}_{X_0}\tilde D=1$. Thus \eqref{p1:cat-superderivative}
and the normalization in Remark~\ref{p1:normalizedpolys} yield
\[
4=2\bigl(\operatorname{ord}_{X_0}\tilde P^2-1\bigr),
\qquad\text{hence}\quad \operatorname{ord}_{X_0}\tilde P^2=3.
\]
Consequently, $(X^d-1)^3$ divides $\tilde P^2$, while $X^d-1$ divides
$\tilde D$ exactly once.
Hence
\begin{align}
g_z^{12}&=\zeta(\lambda_1^2-1)(\lambda_2^2-1)
\left(\frac{\tilde P^2(X)}{\tilde D(X)}\right)^2,\label{p1:cat-normalderivative}\\
&=\zeta(\lambda_1^2-1)(\lambda_2^2-1)(X^d-1)^4
\frac{R(X)^2}{S(X)^2},\qquad X=\frac\zeta\alpha,\label{p1:cat-factorization}
\end{align}
where
\[
R:=\frac{\tilde P^2}{(X^d-1)^3},\qquad
S:=\frac{\tilde D}{X^d-1}
\]
are polynomials in $X$, with $\deg R=\lambda_1+\lambda_2-3d$.

The following lemma distinguishes the roots of $R$, which give planar ends, from those of $S$, where the metric is regular.

\begin{lemma}\label{p1:rootorders}
Assume that~\eqref{p1:cat-normalization} holds.
\begin{enumerate}
\item[(i)] If $X^d=1$, then $R(X)\ne0$ and $S(X)\ne0$.
\item[(ii)] All roots of $S$ are simple and are not roots of $R$. At these points,
which are the poles of $g^{12}$ in $\mathbb C^*$, the metric of $f^{12}$ is regular.
\item[(iii)] All roots of $R$ are simple.
\end{enumerate}
\end{lemma}
\begin{proof}
(i) Let $X_0^d=1$ and write $\mathcal L=X\,d/dX$.
Since $X_0^{\lambda_1}=X_0^{\lambda_2}=1$, the polynomials in Remark~\ref{p1:normalizedpolys} give
\[
\tilde P^2(X_0)=\mathcal L\tilde P^2(X_0)
=\mathcal L^2\tilde P^2(X_0)=0,
\]
\[
\mathcal L^3\tilde P^2(X_0)
=\lambda_1\lambda_2(\lambda_1+\lambda_2)(\lambda_1-\lambda_2)\ne0,
\qquad
\mathcal L\tilde D(X_0)=2\mathcal L^3\tilde P^2(X_0).
\]
Thus $\tilde P^2$ has a zero of order exactly $3$ and
$\tilde D$ a simple zero there. Division by $(X^d-1)^3$ and
$X^d-1$, respectively, proves $R(X_0)\ne0$ and $S(X_0)\ne0$.
In \eqref{p1:cat-normalderivative}, their quotient squared consequently has order $4$, as in
Proposition~\ref{p1:branchorders}(b).

(ii) Let $S(X_0)=0$. By (i), $X_0^d\ne1$, and
\[
\tilde D(X_0)=(X_0^d-1)S(X_0)=0.
\]
By Lemma~\ref{p1:commonroots} and Proposition~\ref{p1:commoncount}, the common roots of $\tilde P^1$ and
$\tilde D$ are exactly the points with $X^d=1$. Since $X_0^d\ne1$, we obtain
$\tilde P^1(X_0)\ne0$. Hence $X_0$ is a pole of $g^{12}$ whose order equals the
multiplicity of $X_0$ as a root of $\tilde D$.
Moreover, $g^1\ne g^2$ at $X_0$:
\begin{itemize}
\item if $g^1=g^2$ were finite at $X_0$, then $g^{12}=g^1$ would be finite there by \eqref{p1:superposition};
\item $g^1=g^2=\infty$ is impossible, since under \eqref{p1:cat-normalization} the poles of $g^i$ lie on the circle
\[
|X|=\left(\frac{\lambda_i-1}{\lambda_i+1}\right)^{1/\lambda_i},
\]
and the function $\lambda\mapsto\lambda^{-1}\log((\lambda-1)/(\lambda+1))$
is strictly increasing on $(1,\infty)$.
\end{itemize}
Hence, by Corollary~\ref{p1:regularcor}, the pole of $g^{12}$ at $X_0$ is simple, i.e. $X_0$ is a simple
root of $\tilde D$. Then $g_z^{12}$ has a pole of order exactly $2$ at $X_0$, and
\eqref{p1:cat-normalderivative} forces
\[
\tilde P^2(X_0)\ne0,\qquad\text{i.e.}\quad R(X_0)\ne0.
\]
In particular, the metric of $f^{12}$ is regular at $X_0$, so $X_0$ is not an end.

(iii) Let $R(X_0)=0$. By (i) and (ii),
\[
X_0^d\ne1,\qquad S(X_0)\ne0,
\]
so $g^{12}$ is finite and branched at $X_0$. By Proposition~\ref{p1:branchorders}, $g^1=g^2$ at $X_0$.
Since the points with $g^1=g^2=g$ are exactly those with $X^d=1$ (Lemma~\ref{p1:commonroots}), we are in
case (a) of Proposition~\ref{p1:branchorders}. Hence $g_z^{12}$ has a zero of order exactly $2$ at $X_0$,
and by \eqref{p1:cat-factorization}
\[
2=\operatorname{ord}_{X_0}g_z^{12}=2\operatorname{ord}_{X_0}R,
\]
so $X_0$ is a simple root of $R$.\end{proof}

If~\eqref{p1:cat-normalization} holds for no $\alpha$, then $g^1=g^2=g$ never occurs.
Thus $\ell=0$ by Lemma~\ref{p1:commonroots}, so $\deg g^{12}=\lambda_1+\lambda_2+1$,
and all branch points in $\mathbb C^*$ are of type~(a) in
Proposition~\ref{p1:branchorders}. They are the roots of $P^2$; for special
$C_{\lambda_i}$ these include common poles of $g^1,g^2$, where $P^2$ and $D$
both vanish.

These roots are simple. Indeed, put
\[
y_1=C_{\lambda_1}\zeta^{\lambda_1},\qquad y_2=C_{\lambda_2}\zeta^{\lambda_2}.
\]
Eliminating $y_1$ from $P^2=\zeta\,dP^2/d\zeta=0$ gives
\[
\lambda_1(\lambda_1^2-1)(\lambda_1-\lambda_2)(\lambda_1+\lambda_2)
\bigl((\lambda_2-1)y_2+\lambda_2+1\bigr)^2=0.
\]
Thus $y_2=-(\lambda_2+1)/(\lambda_2-1)$, and together with $P^2=0$ this forces
$y_1=-(\lambda_1+1)/(\lambda_1-1)$, i.e. $g^1=g^2=g$, which is excluded.

We now translate these root counts into statements about the planar ends of $f^{12}$ and their multiplicities.

\begin{figure}[H]
\centering
\includegraphics[width=0.8\linewidth, height=7cm, keepaspectratio]{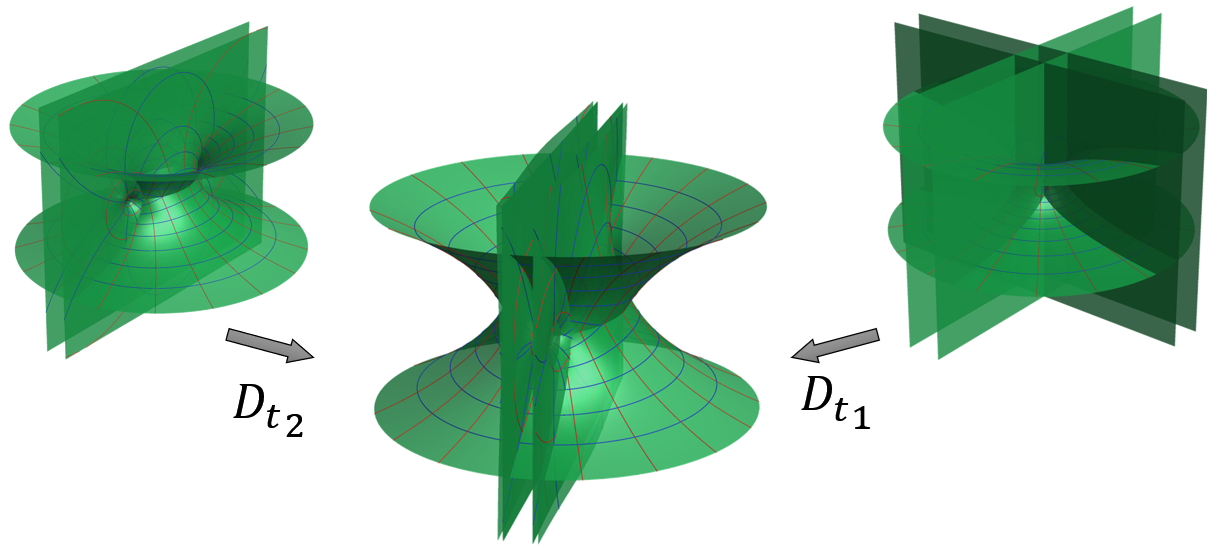}
\caption{Superposition principle for minimal Darboux transforms of the catenoid. Left and right: the single minimal Darboux transforms $f^1=D_{t_1}(f)$ and $f^2=D_{t_2}(f)$, with the parameters chosen so that $\lambda_1=2$ and $\lambda_2=4$. Middle: the superposed surface $f^{12}$ in case~(i) of Proposition~\ref{p1:global}. In particular, since $\lambda_2=2\lambda_1$, all planar ends of $f^{12}$ have multiplicity~3, and there are $d=\gcd(\lambda_1,\lambda_2)=2$ of them.}\label{p4-Catenoid}
\end{figure}

\begin{figure}[H]
\centering
\includegraphics[width=0.8\linewidth]{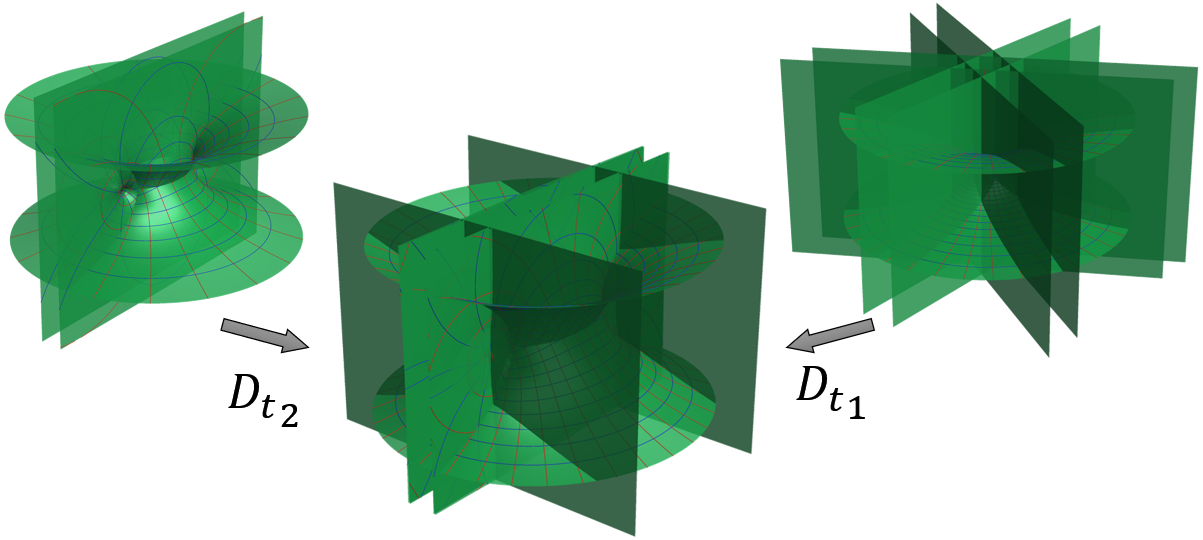}
\caption{Left and right: the single minimal Darboux transforms $f^1=D_{t_1}(f)$ and $f^2=D_{t_2}(f)$, with the parameters chosen so that $\lambda_1=2$ and $\lambda_2=6$. Middle: the superposed surface $f^{12}$ in case~(i) of Proposition~\ref{p1:global}. In this case, $f^{12}$ has $d=\gcd(\lambda_1,\lambda_2)=2$ planar ends of multiplicity~3 and $\lambda_1+\lambda_2-3d=2+6-3\cdot2=2$ embedded planar ends.}\label{p5-Catenoid}
\end{figure}

\begin{figure}[H]
\centering
\includegraphics[width=0.8\linewidth]{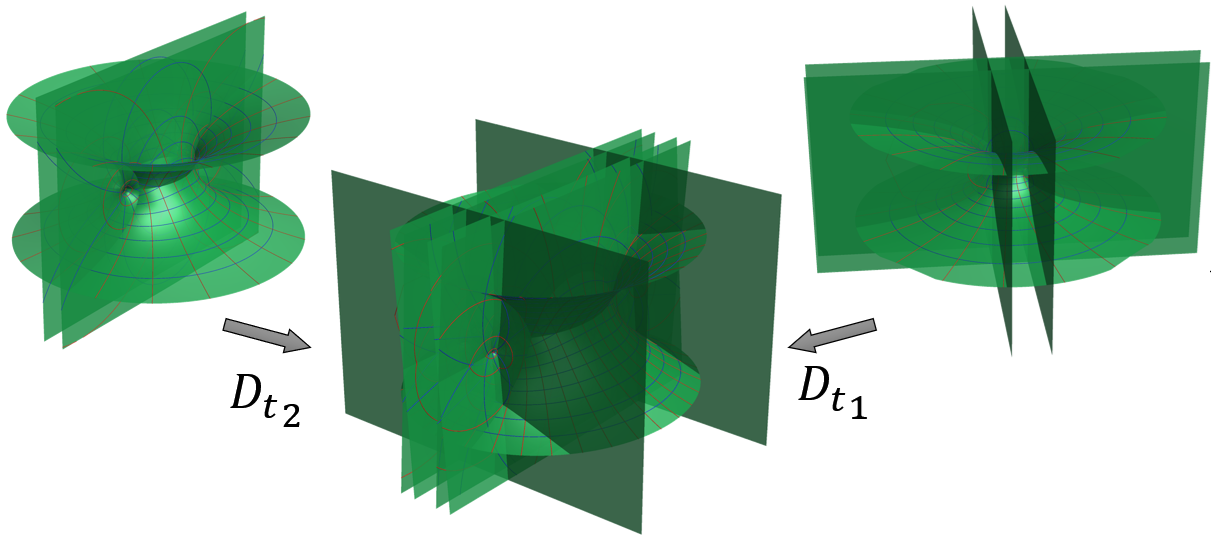}
\caption{Left and right: the single minimal Darboux transforms $f^1=D_{t_1}(f)$ and $f^2=D_{t_2}(f)$, with the parameters chosen so that $\lambda_1=2$ and $\lambda_2=4$. Middle: the superposed surface $f^{12}$ in case~(ii) of Proposition~\ref{p1:global}. In this case, all $\lambda_1+\lambda_2=6$ planar ends of $f^{12}$ are embedded.}
\label{p6-Catenoid}
\end{figure}

\begin{remark}\label{p1:well-defined}
\textup{(a) Well-definedness.} Proposition~\ref{p1:cdprop} gives
\[
f^i=f+\frac1{t_i}(f^{i*}-f^*)^{-1},\qquad
f^{12}=f^1+\frac1{t_2}(f^{12*}-f^{1*})^{-1}.
\]
These expressions are algebraic in $2\pi\mathrm i$-periodic data, so the surfaces are single-valued on $\mathbb C^*$ with their ends removed.

\textup{(b) Planar ends.} By~(a) and Corollary~\ref{p1:endcor}, condition~\eqref{p1:cat-normalization} gives embedded planar ends at the roots of $R$ and ends of multiplicity~3 at the $d$ points $X^d=1$. If~\eqref{p1:cat-normalization} holds for no $\alpha$, the $\lambda_1+\lambda_2$ roots of $P^2$ give embedded planar ends.

\textup{(c) Consistency check.} With two catenoidal ends and $\chi(\overline M)=2$, the Jorge--Meeks formula gives
\[
\int K\,dA=-4\pi(\lambda_1+\lambda_2+1-d)
\]
from $\lambda_1+\lambda_2-3d$ embedded planar ends and $d$ planar ends of multiplicity~3 (see Figures~\ref{p4-Catenoid} and \ref{p5-Catenoid}). If~\eqref{p1:cat-normalization} holds for no $\alpha$, it gives $-4\pi(\lambda_1+\lambda_2+1)$ (see Figure~\ref{p6-Catenoid}). Both values agree with the Gauss-map degree calculation.
\end{remark}

To obtain the global classification, we examine the punctures at $\zeta=0$ and $\infty$ and prove completeness.

\begin{proposition}\label{p1:global}
Let $\lambda_1\ne\lambda_2$ be integers at least~2, let
$C_{\lambda_1},C_{\lambda_2}\in\mathbb C^*$, and set
\[
M_{\lambda_1,\lambda_2}:=\mathbb{CP}^1\setminus
\bigl(\{0,\infty\}\cup\{\zeta\in\mathbb C^*:P^2(\zeta)=0\}\bigr).
\]
The minimal surface $f^{12}$ on $M_{\lambda_1,\lambda_2}$ with Gauss map
\eqref{p1:cat-super} and Weierstrass data $(g^{12},-dz/g_z^{12})$ is complete,
has finite total curvature, and has two catenoidal ends at $\zeta=0,\infty$.
Moreover:
\begin{enumerate}
\item[(i)] If~\eqref{p1:cat-normalization} holds, there are $d=\gcd(\lambda_1,\lambda_2)$
planar ends of multiplicity~3 at $\zeta^d=\alpha^d$ and
$\lambda_1+\lambda_2-3d$ embedded planar ends at the roots of $R(\zeta/\alpha)$
(none if $\{\lambda_1,\lambda_2\}=\{d,2d\}$). The total curvature is
\[
\int_{M_{\lambda_1,\lambda_2}}K\,dA=-4\pi(\lambda_1+\lambda_2+1-d),
\]
and the total number of ends is $\lambda_1+\lambda_2-2d+2$.
\item[(ii)] Otherwise, there are $\lambda_1+\lambda_2$ embedded planar ends at
the roots of $P^2$. The total curvature is $-4\pi(\lambda_1+\lambda_2+1)$,
and the total number of ends is $\lambda_1+\lambda_2+2$.
\end{enumerate}
\end{proposition}

\begin{figure}[H]
\centering
\includegraphics[width=\linewidth]{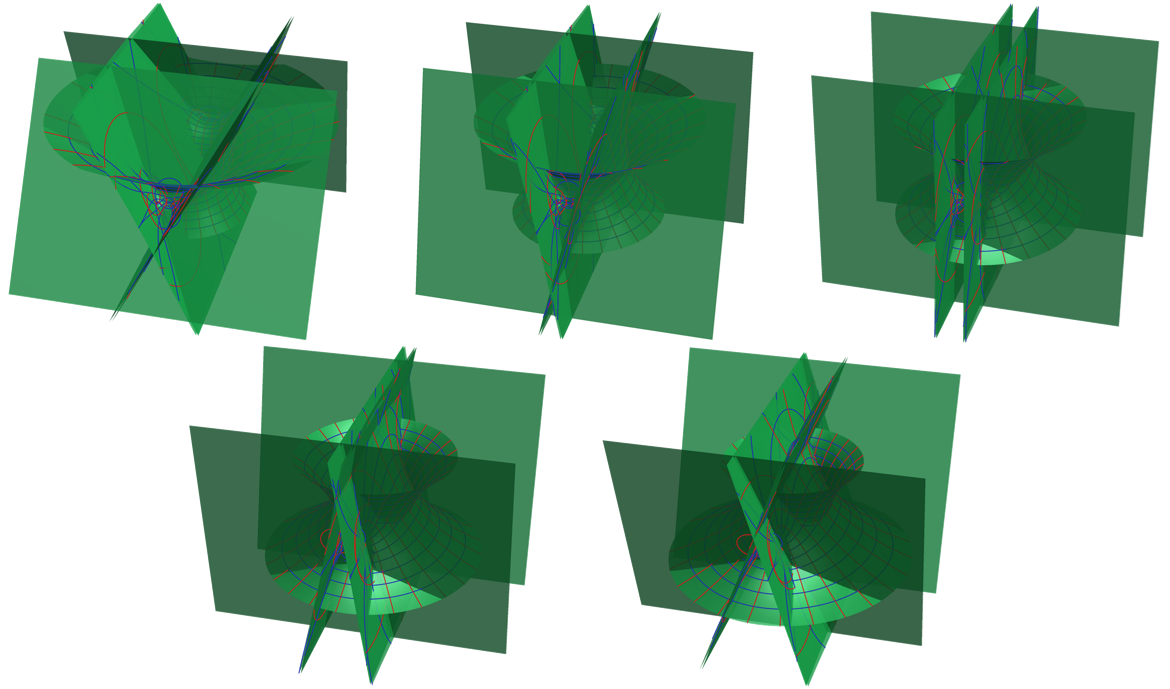}
\caption{The L\'opez--Ros deformation $\{f^{12}_r\}_{r>0}$ of the superposed surface $f^{12}_1$ with $\lambda_1=2$ and $\lambda_2=6$ in case~(i) of Proposition~\ref{p1:global} (Proposition~\ref{p1:lopezros}). From the top left, row by row from left to right, $r=0.6,\,0.8,\,1,\,1.2,\,1.4$.}\label{p7-Catenoid}
\end{figure}

\begin{proof}
The surface is single-valued by Remark~\ref{p1:well-defined}(a).
Under~\eqref{p1:cat-normalization},
$P^2(\zeta)=(\lambda_1+1)(\lambda_2+1)(X^d-1)^3R(X)$, where $X=\zeta/\alpha$.
\begin{itemize}
\item \textit{Regularity.} In case~(i), Proposition~\ref{p1:branchorders} and
Lemma~\ref{p1:rootorders} identify the branch points in $\mathbb C^*$ with
the roots of $P^2$; all poles of $g^{12}$ are simple. In case~(ii), the
calculation before Remark~\ref{p1:well-defined} identifies the branch points with the roots of $P^2$.
At every other pole $g^1\ne g^2$, so Corollary~\ref{p1:regularcor} gives
a simple pole. Hence $(g^{12})^\sharp\ne0$ on $M_{\lambda_1,\lambda_2}$.

\item \textit{Catenoidal ends.} The constant and leading coefficients of
$P^1,D$ give, with $\lambda_*:=\min(\lambda_1,\lambda_2)\ge2$,
\[
g^{12}(\zeta)=
\begin{cases}
a_0\zeta\bigl(1+O(\zeta^{\lambda_*})\bigr),&\zeta\to0,\\
a_\infty\zeta\bigl(1+O(\zeta^{-\lambda_*})\bigr),&\zeta\to\infty,
\end{cases}
\quad
a_0=\frac{(\lambda_1+1)(\lambda_2+1)}
{(\lambda_1-1)(\lambda_2-1)},\quad a_\infty=a_0^{-1}.
\]
Write $\eta=h\,dz=-d\zeta/(\zeta^2\partial_\zeta g^{12})$.
At $0$ these expansions imply
\[
\eta=-\frac{1+O(\zeta^{\lambda_*})}{a_0\zeta^2}\,d\zeta,
\qquad
g^{12}\eta=-\frac{d\zeta}{\zeta}
+O(\zeta^{\lambda_*-1})\,d\zeta.
\]
Thus the horizontal primitive has a simple pole and
$f_3=-\log|\zeta|+O(1)$. With $w=1/\zeta$ at the other end,
$g^{12}=a_\infty w^{-1}(1+O(w^{\lambda_*}))$ and
$\eta=a_\infty^{-1}(1+O(w^{\lambda_*}))\,dw$ give the same
catenoidal behavior. This also verifies directly the conclusion of
\citep[Theorem 3.7]{MRT}.
\item \textit{Planar ends.} Apply Remark~\ref{p1:well-defined}(b), namely
Corollary~\ref{p1:endcor} together with Lemma~\ref{p1:rootorders} in case~(i)
and the preceding root calculation in case~(ii).

\item \textit{Completeness and total curvature.} The line element is
$ds=(1+|g^{12}|^2)|\eta|/2$.
The preceding expansions give $ds\asymp|\zeta|^{-2}|d\zeta|$ at $0$
and $ds\asymp|w|^{-2}|dw|$ at $\infty$; Proposition~\ref{p1:planarprop} treats the
planar ends. Since there are no other singularities, the surface is complete.
For $\zeta=x+\mathrm i y$, the Gauss map gives
\[
-K\,dA=\frac{4|\partial_\zeta g^{12}|^2}{(1+|g^{12}|^2)^2}\,dx\,dy.
\]
A rational map covers the sphere $\deg g^{12}$ times, counting
multiplicity. Its spherical area is therefore $4\pi\deg g^{12}$,
where the degree is $\lambda_1+\lambda_2+1-d$ in case~(i) and
$\lambda_1+\lambda_2+1$ in case~(ii). This proves finite total curvature and
the asserted formulas without presupposing the hypotheses of Osserman's theorem.
\end{itemize}
\end{proof}
Under~\eqref{p1:cat-normalization}, the polynomials $\tilde P^1,\tilde D$ and $R$ are independent of $\alpha$. Set
\[
\Gamma(X):=\frac{X\tilde P^1(X)}{\tilde D(X)},
\qquad g^{12}(\zeta)=\alpha\Gamma(\zeta/\alpha).
\]
Let $f_\alpha^{12}$ be the surface of Proposition~\ref{p1:global}(i), determined up to translation. In the coordinate $X=\zeta/\alpha$, its domain
\[
M=\mathbb{CP}^1\setminus\bigl(\{0,\infty\}\cup\{X:X^d=1\}\cup\{X:R(X)=0\}\bigr)
\]
is independent of $\alpha$. Since $g_z^{12}=\alpha X\Gamma'(X)$ and $dz=dX/X$, its Weierstrass data are
\begin{equation}
(g,\eta)=\bigl(\alpha\Gamma(X),\alpha^{-1}\tilde h(X)\,dX\bigr),
\qquad\tilde h(X):=-\frac1{X^2\Gamma'(X)}.\label{p1:cat-lopezros}
\end{equation}
Here $\eta=h\,dz$ denotes the Weierstrass one-form, independently of the coordinate. Comparing these data with those at $\alpha=1$ separates the geometric effects of the modulus and argument of $\alpha$ (see Figure~\ref{p7-Catenoid}).

\begin{proposition}\label{p1:lopezros}
Write $\alpha=re^{\mathrm i\vartheta}$ with $r>0$. After compatible choices of translation, $f_{\alpha}^{12}=R_{\vartheta} f_r^{12}$, where $R_{\vartheta}$ is the rotation through $\vartheta$ about the $x_3$-axis.
\end{proposition}

\begin{proof}
By \eqref{p1:cat-lopezros}, the data of $f_\alpha^{12}$ are
$(e^{\mathrm i\vartheta}r\Gamma,e^{-\mathrm i\vartheta}r^{-1}\tilde h\,dX)$.
The substitution $(g,h)\mapsto(e^{\mathrm i\vartheta}g,e^{-\mathrm i\vartheta}h)$ is
the rotation through $\vartheta$ about the $x_3$-axis, and $(g,h)\mapsto(rg,r^{-1}h)$
with $r>0$ is the L\'opez--Ros deformation.\end{proof}

\section*{Acknowledgements}
The authors thank Shintaro Akamine, Joseph Cho, Shoichi Fujimori and Wayne Rossman for valuable comments on an earlier version of this work. They are also grateful to the audience at the 3rd OCAMI-IMAG Joint Conference on Differential Geometry for stimulating discussions.

\section*{Statements and declarations}
\subsection*{Funding}
The second author was partially supported by JSPS KAKENHI Grant Number JP25K06978 and by the Osaka Central Advanced Mathematical Institute
(MEXT Promotion of Distinctive Joint Research Center Program
JPMXP0723833165), Osaka Metropolitan University.
\subsection*{Competing interests}
The authors have no relevant financial or non-financial interests to disclose.
\subsection*{Author contributions}
Both authors contributed equally to this work.
\subsection*{Data availability}
Data sharing is not applicable to this article as no datasets were generated or analyzed during the current study.
\subsection*{Use of generative artificial intelligence}
During the preparation of this work, the large language models ChatGPT (OpenAI) and Claude (Anthropic) were used to assist in preparing the exposition, checking computations, investigating observations arising from the examples and developing proofs of several properties. All calculations and arguments were independently checked. The authors take full responsibility for the mathematical content of the article.

\bigskip
\noindent\begin{minipage}{\linewidth}
\small
\textsc{Sakuma Takeshita}\\
Graduate School of Sciences and Technology for Innovation\\
Tokushima University\\
2-1 Minamijosanjima, Tokushima 7708506, Japan\\
\emph{Email address:} \href{mailto:c612537008@tokushima-u.ac.jp}{\texttt{c612537008@tokushima-u.ac.jp}}
\end{minipage}

\medskip
\noindent\begin{minipage}{\linewidth}
\small
\textsc{Masashi Yasumoto}\\
Graduate School of Technology, Industrial and Social Sciences\\
Tokushima University\\
2-1 Minamijosanjima, Tokushima 7708506, Japan\\
\emph{Email address:} \href{mailto:yasumoto.masashi@tokushima-u.ac.jp}{\texttt{yasumoto.masashi@tokushima-u.ac.jp}}
\end{minipage}

\end{document}